%% file: CombCNPC.tex
\documentclass[ 11pt]{article}

\usepackage[utf8]{inputenc}  
\usepackage[T1]{fontenc}         % Police contenant les caractÃƒÆ’Ã†â€™Ãƒâ€šÃ‚Â¨res franÃƒÆ’Ã†â€™Ãƒâ€šÃ‚Â§ais
\usepackage{geometry}         % DÃƒÆ’Ã†â€™Ãƒâ€šÃ‚Â©finir les marges
\usepackage[english]{babel}  % Placez ici une liste de langues, la
\usepackage{cite}
\usepackage{hyperref}
\usepackage{float}
\usepackage{amssymb}   
\usepackage[titletoc,title]{appendix}
\usepackage{amsmath} 
\usepackage{amscd}
\usepackage{graphicx, color}    
\usepackage{amsthm}                 
\usepackage{latexsym}     % derniÃƒÆ’Ã†â€™Ãƒâ€šÃ‚Â¨re ÃƒÆ’Ã†â€™Ãƒâ€šÃ‚Â©tant la langue principale

\usepackage[all]{xy}
\usepackage{verbatim}
\usepackage{enumerate}
\newtheorem{thm}{Theorem}[section]
\newtheorem*{thm*}{Theorem}
\newtheorem*{thmA}{Theorem A}
\newtheorem*{thmB}{Theorem B}
\newtheorem*{thmC}{Theorem C}
\newtheorem*{thmD}{Theorem D}

\newtheorem{cor}[thm]{Corollary}
\newtheorem*{cor*}{Corollary}
\newtheorem*{question}{Question}
\newtheorem{prop}[thm]{Proposition}
\newtheorem{lem}[thm]{Lemma}% \pagestyle{headings}        % Pour mettre des entÃƒÆ’Ã†â€™Ãƒâ€šÃ‚Âªtes avec les titres

\theoremstyle{definition}
\newtheorem{definition}[thm]{Definition}
\newtheorem{deflem}[thm]{Lemma-Definition}

\newtheorem{rmk}[thm]{Remark}
\newtheorem{conv}[thm]{Convention}
\newtheorem{exm}[thm]{Example}
\newtheorem{notation}[thm]{Notation}

\def\lquotient#1#2{%
\makeatletter
\lower.6ex\hbox{$#1$}\backslash\raise.3ex\hbox{$#2$}%
\makeatother
}															

\def\rquotient#1#2{%
\makeatletter
\raise.6ex\hbox{$#1$}/\raise.3ex\hbox{$#2$}%
\makeatother
}			
\newcommand{\CG}{{\cC G}}
\newcommand{\VG}{{\cV G}}

\newcommand{\cB}{{\mathcal B}}
\newcommand{\cC}{{\mathcal C}}

\newcommand{\cE}{{\mathcal E}}

\newcommand{\cP}{{\mathcal P}}

\newcommand{\cR}{{\mathcal R}}

\newcommand{\cU}{{\mathcal U}}
\newcommand{\cV}{{\mathcal V}}

\newcommand{\cY}{{\mathcal Y}}
\newcommand{\cZ}{{\mathcal Z}}

\newcommand{\ra}{\rightarrow}

\newcommand{\hra}{\hookrightarrow}

\newcommand{\bd}{\partial}
\newcommand{\bds}{\partial_{Stab} G}

\newcommand\blfootnote[1]{%
	\begingroup
	\renewcommand\thefootnote{}\footnote{#1}%
	\addtocounter{footnote}{-1}%
	\endgroup
}

\newcommand{\wt}{\widetilde}

\title{\textbf{A combination theorem for combinatorially non-positively curved complexes of  hyperbolic groups}}           
\author{Alexandre Martin and Damian Osajda}
\date{}
\newcommand{\Addresses}{{% additional braces for segregating \footnotesize
		\bigskip
		\footnotesize
		
		Alexandre Martin, \textsc{Department of Mathematics, Heriot-Watt University
			EH14 4AS Edinburgh, UK}\par\nopagebreak
		\texttt{alexandre.martin@hw.ac.uk}
		
		\medskip
		
		Damian Osajda, \textsc{Instytut Matematyczny,  UWr., pl.\ Grunwaldzki 2/4, 
			50--384 Wroc\-{\l}aw, Poland}\par\nopagebreak
		\textsc{{Institute of Mathematics, Polish Academy of Sciences, \'Sniadeckich 8,  
				00--656 War\-sza\-wa, Poland}}\par\nopagebreak
		\texttt{dosaj@math.uni.wroc.pl}
		
}}

\begin{document}
\maketitle

\begin{abstract} 
We prove a combination theorem for hyperbolic groups, in the case of groups acting on  complexes displaying combinatorial features reminiscent of non-positive curvature. Such complexes include for instance weakly systolic complexes and $C'(1/6)$ small cancellation polygonal complexes. Our proof involves constructing 
a potential Gromov boundary for the resulting groups and analyzing the dynamics of the action on the boundary
in order to use Bowditch's characterization of hyperbolicity. A key ingredient is the introduction of a combinatorial property that implies a weak form of non-positive curvature, and which  holds for large classes of complexes. 

 As an application, we study the hyperbolicity of groups obtained by small cancellation over a graph of hyperbolic groups. 
\end{abstract}

%\tableofcontents

\input{intro.tex}

\input{comb.tex}

\input{bdry.tex}

\input{wsys.tex}

\input{sc.tex}

\input{example.tex}

\input{appendices.tex}

\bibliographystyle{plain}
\bibliography{CombinationSystolic}

\Addresses

\end{document}

%% file: intro.tex
\section{Introduction}
\label{s:intro}

\emph{Hyperbolic groups}, known also as \emph{Gromov hyperbolic} or \emph{$\delta$--hyperbolic groups}, were introduced by Gromov \cite{GromovHyperbolicGroups}. This concept unified
approaches to various classical groups that had been studied before. Examples include: free groups, lattices in automorphisms groups of Lobachevski hyperbolic space (called sometimes also `hyperbolic groups'), and many small cancellation groups. 
Over the last thirty years the theory of hyperbolic groups has been at the centre of group theory,
three-manifolds theory, and influenced many other disciplines, including ones  outside mathematics such as 
 computer science.
Besides constructing arithmetic lattices, another most important way of creating hyperbolic groups
is the use of different gluing techniques.
The Bestvina-Feighn combination theorem is considered as the first result of this type \cite{BestvinaFeighnHyperbolicity}. Roughly, it states that some finite graphs of hyperbolic groups
have hyperbolic fundamental groups. Januszkiewicz-{\' S}wi{\k a}tkowski \cite{SystolicComplexes} developed a technique
of constructing high-dimensional hyperbolic groups as fundamental groups of finite complexes of finite groups
satisfying some combinatorial non-positive curvature conditions, called \emph{systolicity}.
Recently, the first author presented an approach unifying in a way the two above constructions  \cite{MartinBoundaries}. He showed that, under some assumptions, a finite non-positively
curved (that is, locally CAT(0)) complex of hyperbolic groups has a hyperbolic fundamental group. In other words, a group acting
co-compactly on a non-positively curved complex with hyperbolic stabilizers and satisfying some acylindricity condition is hyperbolic. 

In the current paper we present a `combinatorial counterpart' of \cite{MartinBoundaries}. The motivation is as 
follows. Non-positive curvature (NPC), in the sense of the local CAT(0) property, is a metric feature that 
is quite difficult to verify in general. Only in some simple instances, e.g.\ CAT(0) cube complexes,
can a standard piecewise Euclidean or hyperbolic (that is coming from $\mathbb{H}^n$) metric  be
shown to satisfy the NPC conditions. Even then, one generally uses some equivalent combinatorial criteria to show this condition.
In many other cases it is not at all clear what would be a candidate for a reasonable CAT(0) metric on a complex.
In the approach of Januszkiewicz-{\' S}wi{\k a}tkowski \cite{SystolicComplexes},  instead of a metric setting one relies on a combinatorial
notion of `non-positive curvature'. The NPC condition is replaced by a simple -- and easily checkable -- local combinatorial condition. 
This is the reason why one can relatively easily construct new interesting
examples of hyperbolic complexes and groups acting on them. For the same reason it is
worth exploring actions on combinatorially non-positively curved complexes with hyperbolic    
stabilizers.  We show that such settings lead to new constructions of hyperbolic groups.
In a way our approach is a natural generalization of the ones of Bestvina-Feighn \cite{BestvinaFeighnCombinationHyperbolic} and
of Januszkiewicz-{\' S}wi{\k a}tkowski \cite{SystolicComplexes}.

We now present the most general result of this paper, which is later tailored to some more specific 
situations. The main technical condition therein -- the \emph{Small Angle Property} (see Definition \ref{def:small_angle_property}) -- can be seen as a weak form
of (combinatorial) non-positive curvature, and mimics in combinatorial settings the behaviour of geodesics in a CAT(0) space. A weakening of this property was introduced by the first author in \cite{MartinTame} and has been used to show the hyperbolic features of several groups: groups of birational automorphisms \cite{MartinTame}, certain Artin groups \cite{MartinChatterjiCriterion}, etc. The Small Angle Property 
is satisfied in many interesting 
cases, including the ones we explore afterwards (see Theorem B and Theorem C below).

\begin{thmA}\label{thm:main}
Let $G$ be a group acting cocompactly, without inversions on a hyperbolic complex $X$ 
with finite intervals and satisfying the Small Angle Property. Assume that the following local conditions are satisfied:
\begin{description}
\item[(L1)]
for every face $\sigma$ of $X$ the stabilizer $G_{\sigma}$ of $\sigma$ is hyperbolic;
\item[(L2)]
for any two faces $\sigma \subseteq \sigma'$ the inclusion $G_{\sigma'} \hookrightarrow G_{\sigma}$ is a quasi-convex
embedding.
\end{description}
Furthermore, assume that the following global conditions are satisfied:
\begin{description}
\item[(G1)]
the action of $G$ on $X$ is weakly acylindrical, that is, there exists an upper bound on the distance of two vertices 
that are both fixed by an infinite subgroup of $G$;
\item[(G2)]
loops in fixed-point sets are contractible, that is, for every subgroup $H$ of $G$, every loop contained in the associated fixed-point set $X^H$ is nullhomotopic.
\end{description}
Then $G$ is hyperbolic and the inclusions $G_{\sigma} \hra G$ are quasi-convex embeddings, for every face $\sigma$. 
\end{thmA}

In Subsection~\ref{s:general} we point out an even more general setting in which an analogous result holds.
In order to avoid dealing with overly technical notations and proofs though, we decided to concentrate on the setting mentioned in Theorem A. We also give in Example \ref{example_noSAP} an instance of an action of a non-hyperbolic group on a hyperbolic complex without the Small Angle Property satisfying conditions (L1), (L2), (G1), and (G2). This shows how crucial having a control on the geodesics of the complex is to obtain such combination theorems. 

The approach followed to prove Theorem A is a dynamical one, and goes back to work of Dahmani on graphs of relatively hyperbolic groups \cite{DahmaniCombination}, work of the first author on CAT(0) complexes of hyperbolic groups \cite{MartinBoundaries}, and recent work on CAT(0) complexes of relatively hyperbolic groups by Pal--Paul \cite{PalPaul}. In a nutshell, we construct a candidate $\partial G$ for the Gromov boundary of the group $G$, by gluing together the Gromov boundary of $X$ and the Gromov boundaries of the various stabilisers of vertices, and we endow this set with an appropriate topology. We then study the dynamics of the action of $G$ on $\partial G$, and show that $G$ acts as a uniform convergence group on it (see Section \ref{s:conv}), which implies the hyperbolicity of $G$ and that $\partial G$ is equivariantly homeomorphic to the Gromov boundary of $G$, by a characterisation due to Bowditch \cite{BowditchTopologicalCharacterization}. The construction of the compactification and its topology are similar to the constructions  in \cite{MartinBoundaries}, and the heart of the article is to construct appropriate combinatorial analogues of the tools developed therein to study its topology and the dynamics of the action.

We now consider applications of the main theorem above in the case of particular complexes.
As noted above, systolicity is a well-known instance of a combinatorial non-positive curvature.
In \cite{O-cnpc, ChO} the notion of \emph{weakly systolic} complexes was introduced.
The definition is very close to systolicity, but the class of resulting complexes is very different.
Let us just note here that systolic complexes exhibit some asymptotically `two-dimensional' behaviour -- 
at large scale they do not contain spheres. Such restrictions do not exist for weakly systolic complexes,
although the methods used for exploring both classes are very similar. We show that weakly systolic complexes 
have \emph{tight hexagons} (see Subsection~\ref{s:wsys}), which implies the Small Angle Property, and therefore we obtain the following theorem, which may be seen as a straightforward generalization of Januszkiewicz-{\' S}wi{\k a}tkowski constructions from \cite{SystolicComplexes} -- the finite groups being replaced by general hyperbolic groups.

\begin{thmB}\label{t:wsys}
Let $G$ be a group acting cocompactly, without inversions on a hyperbolic weakly systolic complex $X$ without infinite simplices. Suppose that the conditions (L1), (L2), and (G1) of Theorem A are satisfied. 
%Assume that the following conditions are satisfied:
%\begin{itemize}
%	\item
%	for every simplex $\sigma$ of $X$ the stabilizer $G_{\sigma}$ of $\sigma$ is hyperbolic;
%	\item
%	for any two simplices $\sigma \subseteq \sigma'$ the inclusion $G_{\sigma'} \hookrightarrow G_{\sigma}$ is a quasi-convex
%	embedding.
%	\item
%	the action of $G$ on $X$ is acylindrical, that is, there exists an upper bound on the distance of two vertices 
%	that are both fixed by an infinite subgroup of $G$;
%\end{itemize}
Then $G$ is hyperbolic and the inclusions $G_{\sigma} \hra G$ are quasi-convex embeddings.
\end{thmB}

%The theorem above may be seen as a straightforward generalization of Januszkiewicz-{\' S}wi{\c a}tkowski constructions from \cite{SystolicComplexes} -- the finite groups are replaced by more general hyperbolic groups.

Another important class of combinatorially non-positively curved complexes is the class of \emph{small cancellation complexes}. In this article we focus on the metric version of small cancellation -- \emph{$C'(1/6)$ small cancellation complexes}, called simply \emph{$C'(1/6)$ complexes} (see Subsection~\ref{s:smallc} for precise definitions). Small cancellation complexes and groups
are among the most classical examples of hyperbolic spaces and groups.

\begin{thmC}\label{t:sc}	
	Let $G$ be a group acting cocompactly, without inversions on a $C'(1/6)$ complex $X$, 
	so that the conditions (L1), (L2), and (G1) of Theorem A are satisfied. 
%	Assume that the following conditions are satisfied:
%	\begin{itemize}
%		\item
%		for every face $\sigma$ of $X$ the stabilizer $G_{\sigma}$ of $\sigma$ is hyperbolic;
%		\item
%		for any two faces $\sigma \subseteq \sigma'$ the inclusion $G_{\sigma'} \hookrightarrow G_{\sigma}$ is a quasi-convex
%		embedding.
%		\item
%		the action of $G$ on $X$ is acylindrical, that is, there exists an upper bound on the distance of two vertices 
%		that are both fixed by an infinite subgroup of $G$;
%	\end{itemize}
	Then $G$ is hyperbolic and the inclusions $G_{\sigma} \hra G$ are quasi-convex embeddings, for every face $\sigma$.  
\end{thmC}

Again, the above result may be seen as a generalization of the classical fact, that groups acting geometrically
on $C'(1/6)$ complexes are hyperbolic.

Finally, we apply Theorem C in a specific situation of small cancellation over graphs of groups. The hyperbolicity of certain small cancellation groups over graphs of groups was already considered by the first author in \cite{MartinSmallCancellationClassifying}. However, the small cancellation condition used therein was much stronger, in order to endow some of the spaces considered with a CAT(0) metric. In particular, this stronger condition, generally referred to as $C''(1/6)$, prevents the construction of infinitely presented small cancellation groups in the classical setting. By contrast, the combinatorial approach used here allows us to work with the combinatorial geometry of $C'(1/6)$ polygonal complexes, a much weaker small cancellation condition. While we focus here on quotients obtained by taking finitely many relations, the approach followed in this paper can thus be seen as paving the way for a geometric study of infinitely presented small cancellation groups over graphs of hyperbolic groups.

\begin{thmD}\label{thm:hyperbolicity_small_cancellation}
	Let $G(\Gamma)$ be a finite graph of groups over a simplicial graph $\Gamma$ satisfying the following: 
	\begin{itemize}
		\item every vertex group is hyperbolic,
		\item every edge group embeds as a quasi-convex subgroup in the associated vertex groups,
		\item for every vertex $v$ of $\Gamma$, the family of adjacent edge groups $(G_e)_{v \in e}$ is almost malnormal.
	\end{itemize}
	Let $G$ be the fundamental group of $G(\Gamma)$. Let $\cR$ be a finite set of relators satisfying the classical $C'(1/6)$--small cancellation over $G(\Gamma)$. Then the quotient group $\rquotient{G}{\ll\cR\gg }$ is hyperbolic and the quotient map $G \ra \rquotient{G}{\ll\cR\gg }$  embeds every local group of $G(\Gamma)$ as a quasi-convex subgroup.  
\end{thmD}

We believe that the crucial features used in the proof -- particularly the Small Angle Property
%, having  fixed-point sets with contractible loops, and finite intervals -- 
-- hold for many other classes of  `combinatorially non-positively curved' complexes, and hence similar results for corresponding complexes of groups can be proved along the same lines
using Theorem A. An example of a very general class of graphs containing $1$--skeleta of weakly systolic complexes and of CAT(0) cubical complexes (for which the corresponding combination theorem holds by \cite{MartinBoundaries})
is the class of weakly modular graphs extensively studied in metric graph theory \cite{BaCh-survey}.
Triangle-square complexes associated to weakly modular graphs have been introduced in \cite{CCHO} and shown
to have numerous non-positive-curvature-like features.  
In particular, we are naturally led to the following question: 

\begin{question}
	Do triangle-square complexes associated to weakly modular graphs have the Small Angle Property?
	\end{question}

{\bf Organization of the paper.} In  Section~\ref{s:comb} we present combinatorial preliminaries for our work. First (Subsection~\ref{s:hyper}), we recall some basic facts about hyperbolic complexes
and group actions, then (Subsection~\ref{s:sap}) we define the main conditions on complexes needed in our approach -- the Small Angle Property (Definition~\ref{def:small_angle_property}) and the property of having tight hexagons (Definition~\ref{d:hex}). Section~\ref{s:constr} is devoted to proving the main Theorem A: In Subsection~\ref{s:constrbdry} we define the Gromov boundary of the ambient group, in Subsection~\ref{s:topology} we define the topology on the boundary, and finally, in Subsection~\ref{s:conv} we explain how the dynamics of the action is used to prove  Theorem A. The proofs of the main results are postponed to an appendix, being natural generalisations of the proofs of \cite{MartinBoundaries}.
% by showing that the group acts on the boundary as a uniform convergence group and
%using Bowditch's characterization.
Section~\ref{s:examples} is devoted to proving Theorem B (Subsection~\ref{s:wsys}), Theorem C (Subsection~\ref{s:smallc}), and Theorem D (Subsection~\ref{s:scog}). Finally, in Appendix A we provide proofs
of results used in Section~\ref{s:constr}. Because the proofs are the natural combinatorial counterparts of the  original proofs from
\cite{MartinBoundaries},  we follow the same structure as much as possible, indicating the issues that have to be adapted.
\medskip

\noindent
{\bf Acknowledgments.} Alexandre Martin was supported by  the ERC grant no. 259527 of G. Arzhantseva, and by the FWF grant M1810-N25. Damian Osajda was supported by (Polish) Narodowe Centrum Nauki, grant no.\ UMO-2015/\-18/\-M/\-ST1/\-00050. This work was partially supported by the grant 346300 for IMPAN from the Simons Foundation and the matching 2015-2019 Polish MNiSW fund.

%% file: comb.tex
\section{Combinatorial geometry of hyperbolic complexes}
\label{s:comb}

\subsection{Hyperbolic complexes and groups acting on them}
We recall here a few definitions and conventions that will be used throughout the article.\\
\label{s:hyper}

\label{s:prel}
\textbf{Complexes.} In this article by a \emph{complex} we mean a CW-complex with polyhedral cells and such that the attaching map of every closed cell is an embedding.
% are 
%\emph{combinatorial}, that is, restricted 
%to the interiors of cells are homeomorphisms onto their images; see e.g.\ \cite{Hat}.
Particular classes of such complexes considered by us will be: simplicial graphs, $C'(1/6)$ polygonal complexes (see Subsection~\ref{s:smallc}), and some
flag simplicial complexes (in Subsection~\ref{s:wsys}). 
For a complex $X$, by $X^{(k)}$ we denote the $k$--skeleton of $X$. We usually assume that $X$, and thus $X^{(2)}$, is connected and simply connected. Throughout this article we further assume that the $1$--skeleton $X^{(1)}$ of any complex $X$ is a \emph{simplicial graph}, that is, a graph without loops and multiple edges. (Observe that this can by assured by subdividing the complex.)

\textbf{Distance and geodesics.} We consider the set $X^{(0)}$ of vertices of $X$ as a metric space with the distance given by the number of edges in the shortest
combinatorial path in $X^{(1)}$ between vertices. Such shortest paths are called \emph{geodesics}. 
Recall that a subcomplex of a complex is \emph{convex} if every geodesic in the $1$--skeleton
joining two vertices of the subcomplex is contained in this subcomplex.
We also consider \emph{oriented geodesics}, that might be thought of as
ordered (in a natural way) sequences of vertices or edges of a geodesic. The \emph{interval}
$I(u,v)$ between two vertices $u$ and $v$ consists of all the vertices on geodesics between $u$ and $v$, that is,  all the vertices $x$ satisfying the equality $d(u,x)+d(x,v)=d(u,v)$.

\textbf{Hyperbolicity.} We say that \emph{$X$ is hyperbolic} whenever $X^{(0)}$ endowed with the above distance is a Gromov hyperbolic metric space. From now on, unless stated otherwise, we will consider only hyperbolic complexes.
By $\partial X$ we denote the Gromov boundary of $X$, and by $\overline X$ we denote
the bordification $X\cup \partial X$. Note that $\overline X$ is metrisable, but is not necessarily compact when $X$ is not locally compact.
By a \emph{generalised vertex of $\overline{X}$} we mean a point of $X^{(0)} \cup \partial X$.

\textbf{Group actions.}  We consider group actions on complexes by cellular automorphisms. Unless stated otherwise, we assume that groups
act \emph{without inversions}, that is, if an element of the group stabilizes a cell set-wise
then it stabilizes the cell point-wise. (Observe that a group acting on a complex acts without inversions on the subdivision of the complex.)
% Alex: This should be a proper definition, it is not satisfied in full generality. We assume that for any subgroup $H<G$ the set of points fixed by $H$ is empty or simply connected.
Of particular interest will be the case where loops in fixed-point sets of the action of $G$ on $X$ are contractible, that is, if for every subgroup $H$ of $G$, every loop contained in the associated fixed-point set $X^H$ is nullhomotopic.
%dam 

%
%Note that, by a slight abuse of notation,  we do not require in the previous definition that fixed-point sets are connected.

\subsection{The Small Angle Property and tight hexagons.}
\label{s:sap}

In this subsection, we introduce a  property which can be thought as a combinatorial counterpart of the properties of geodesics in a CAT(0) space. Namely, we define the Small Angle Property appearing in the formulation of Theorem A, together with another useful property implying it -- having tight hexagons -- that will be used in Section~\ref{s:examples}.

\begin{definition}[exit edge]
Let $K$ be a convex subcomplex of $X$, $v$ a generalised vertex of $\overline{X}$ and $w$ a generalised vertex of $\overline{X}$ which is not in $K$. We say that an oriented geodesic from $v$ to $w$ \emph{goes through} $K$ if it contains a vertex of $K$. 
If an oriented geodesic $\gamma$ goes through $K$, there exists an edge of $\gamma$ not contained in $K$ such that
all consecutive edges are not in $K$ as well. The first such edge is called the \emph{exit edge} and denoted $e_K(\gamma)$.
\end{definition}

\begin{definition}[path around a subcomplex]
Let $K$ be a convex subcomplex of $X$, and $e,e'$ two edges of $X$ such that $e,e' \nsubseteq K$ and $e \cap K, e' \cap K \neq \varnothing$. A \emph{path around $K$} between $e$ and $e'$ is a sequence $\sigma_1, \ldots, \sigma_n$ of cells of $X$ and a sequence $e_0:=e,e_1,\ldots,e_n:=e'$ of edges of $X$ not contained in $K$ such that $e\subset \sigma_1$, $e' \subset \sigma_n$ and for every $i$, $e_i\subset \sigma_i \cap \sigma_{i+1}$, $e_i \cap K \neq \varnothing$. The integer $n$ is called the \emph{length} of this path.
\end{definition}

The following property will be the main ingredient in studying the topology of the boundary $\partial G$ in the Appendix.

\begin{figure}[h]\label{f:sap}
	\begin{center}
		\includegraphics[scale=0.75]{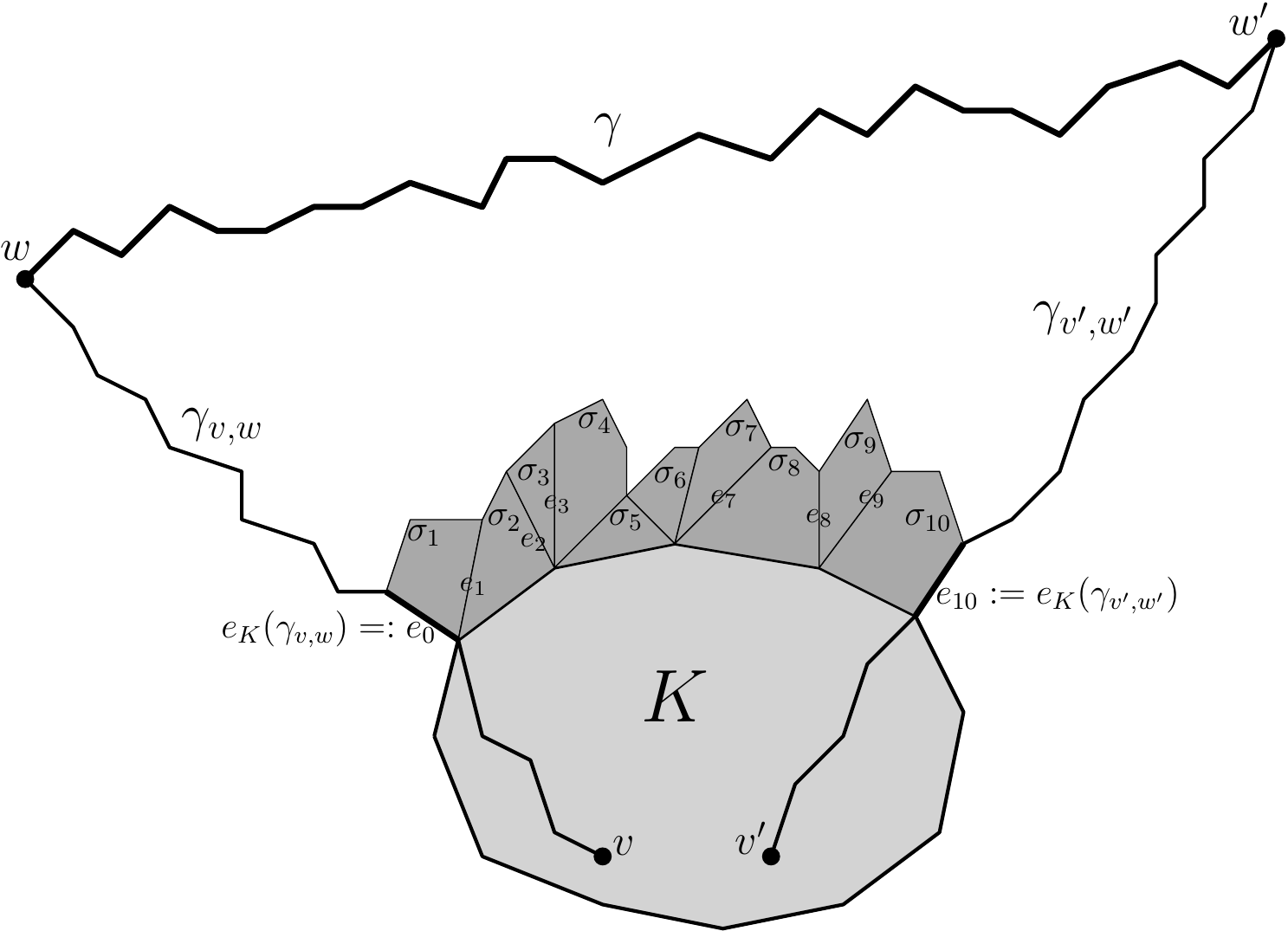}
	\end{center}
	\caption{The Small Angle Property. A path around $K$ between $e_K(\gamma_{v,w})$ and $e_K(\gamma_{v',w'})$ is shaded dark gray.}
\end{figure}

\begin{definition}[Small Angle Property]\label{def:small_angle_property}
The complex $X$ 
is said to satisfy the \textit{Small Angle Property} if for every integer $n \geqslant 0$, there exists an integer $r(n) \geqslant 1$ such that the following holds: 

Let $K$ be a convex subgraph of $X$ with at most $n$ edges. Let  $v,v'$ be vertices of $K$. Let $\gamma$ be a geodesic disjoint from $K$ and let $w,w'$ be generalised vertices of $\gamma$ (that is vertices of $\gamma$ or points in $\partial X$ defined by $\gamma$). Let $\gamma_{v,w}$ and $\gamma_{v',w'}$ be oriented geodesics
from $v$, respectively $v'$, to $w$, respectively $w'$. Then there exists a path around $K$ between $e_K(\gamma_{v,w})$ and $e_K(\gamma_{v',w'})$ of length at most $r(n)$; see Figure~\ref{f:sap}. 
\end{definition}
This notion is a strengthening of the notion of having a \textit{bounded angle of view}, introduced by the first author in \cite{MartinTame}. We now introduce the main combinatorial tool used for proving that a given polyhedral complex satisfies the Small Angle Property.

\begin{definition}[disc diagram, reduced disc diagrams, arcs] A \textit{disc diagram} $D$ of the 
%	$C'(1/6)$ polygonal 
	complex $X$ is a contractible planar polygonal complex endowed with a combinatorial map  $D \ra X$ which is an embedding on each polygon. For a disc diagram $D$, we denote by $\partial D$ its boundary and $\mathring{D}$ its interior. The \textit{area} of a diagram $D$, denoted Area$(D)$, is the number of polygons of $D$. 

For a polygon $R$ of $D$, the intersection $\partial R \cap \partial D$ is called the \textit{outer component} of $R$ (and the \textit{outer path} if such an intersection is connected), the closure of $\partial R \cap \mathring{D}$ is called the \textit{inner component} of $R$ (and the \textit{inner path} if such an intersection is connected).

A diagram is called \textit{non-singular} if its boundary is homeomorphic to a circle, \textit{singular} otherwise. 

A diagram is called \textit{reduced} if no two distinct polygons of $D$ that share an edge are sent to the same polygon of $X$.

The \emph{degree} of a vertex of $D$ is the number of edges containing the vertex. An \textit{arc} of $D$ is a path of $D$ whose interior vertices have degree $2$ and whose boundary vertices have degree at least $3$. Such an arc is called \textit{internal} if its interior is contained in $\mathring{D}$, \textit{external} if  the arc is fully contained in $\partial D$.
\end{definition}

By the relative simplicial approximation theorem \cite{Zeeman1964}, for every cycle in a simply connected complex there exists a disc diagram (cf.\ also van Kampen's lemma e.g.\ in \cite[pp.\ 150-151]{LyndonSchupp}).
%Recall the following fundamental result:
%
%\begin{thm}[Lyndon--van Kampen]
%Every loop of $X$ is the boundary of a reduced disc diagram. \qed
%\end{thm} 

%By a \emph{hexagon} in $X$ we mean a geodesic hexagon, that is, a cycle consisting of 
%six vertices and six geodesics 
%joining them consequently. We allow some of these geodesics to be degenerate, that is, some of the
%vertices may coincide.

\begin{definition}[hexagons, tight hexagons]
\label{d:hex}	
A (geodesic) \emph{hexagon} of $X$ is a cycle consisting of 
six vertices and six geodesics 
joining them consequently. We allow some of these geodesics to be degenerate, that is, some of the
vertices may coincide.

%Let $k \geqslant 3$ be an integer. 
We say that the complex $X$ has \textit{tight hexagons} if there exists an integer $N$ such that every embedded hexagon of $X$ bounds a disc diagram $D$, called a \textit{tight} disc diagram, satisfying the following:
\begin{itemize}
\item the map $D \ra X$ is at most $N$-to-$1$,
\item each vertex of $D$ has degree at most $N$.
\end{itemize}
%Such a disc diagram $D$ is called a \textit{tight} disc diagram.
%% We say that the complex $X$ has  \emph{tight disc diagrams} if it has tight $\cP$-$k$-gons for every $k \geqslant 3$ (with constants $N_k$ that a priori depend on $k$). If this holds, we say that the system of paths $\cP$ is \emph{tight}. 
%We say that the collection of combinatorial geodesics of $X$ is \textit{tight} if $X$ has tight hexagons.
%
%For a polyhedral complex endowed with a tight system of paths $\cP$, we will simply refer to an element of $\cP$ as a \emph{tight path} if no confusion is possible.
\end{definition}

%\begin{rmk}
%It is straightforward to check that having tight $\cP$-triangles implies having tight $\cP$-polygons.
%\end{rmk}

%\begin{rmk}
% ??? Say something about isoperimetric inequality in non-locally finite hyperbolic complexes, with the basic counterexample being the cone on a simplicial line.
%\end{rmk}

%\begin{rmk}
% ??? Say something about tight geodesics in the curve complex and the lemma which is in the spirit of tight discs diagrams.
%\end{rmk}

\begin{lem}
A polyhedral complex $X$ with tight hexagons satisfies the Small Angle Property. \label{tightdiscsimplyrefinement}
\end{lem}

\begin{proof}
Let $K$, $n$, $v,v'$, $\gamma$, $w,w'$ be as in Definition~\ref{def:small_angle_property}. Choose an oriented  geodesic $\gamma_{v,w}$ from $v$ to $w$ (respectively $\gamma_{v',w'}$ from $v'$ to $w'$, respectively $\gamma_{v',v}$ from $v'$ to $v$). Note that $K$ being convex, we have $\gamma_{v,v'} \subset K$. Denote by $\gamma_{w,w'}$ the portion of $\gamma$ between  $w$ and $w'$. If $w$ (respectively $w'$) is a point of the Gromov boundary $\partial X$, we can choose vertices $a, b$ (respectively $a', b'$) of $X$ such that $a \in \gamma_{v,w}, b \in \gamma_{w,w'}$ (respectively $a' \in \gamma_{v',w'}, b' \in \gamma_{w,w'}$) and geodesics $[a,b]$, $[a'b']$ that are disjoint from $K$. 
If $w$ (respectively $w'$) is a vertex of $X$ then we set $a:=w=:b$ (respectively $a':=w'=:b'$).
Let $\gamma_{v,a}$ (respectively $\gamma_{v',a'}$) be the oriented sub-geodesic of $\gamma_{v,w}$ (respectively $\gamma_{v',w'}$) from $v$ to $a$ (respectively $v'$ to $a'$). Let $\gamma_{a,b}$ (respectively $\gamma_{a',b'}$) be an oriented geodesic from $a$ to $b$ (respectively $a'$ to $b'$). Let $\gamma_{b,b'}$ be the oriented sub-geodesic of $\gamma$ from $b$ to $b'$.
Without loss of generality, we can assume that the hexagon $$ \gamma':= \gamma_{v',v} \cup \gamma_{v,a} \cup \gamma_{a,b} \cup \gamma_{b,b'} \cup  \gamma_{b',a'} \cup \gamma_{a',v'} $$ is embedded in $X$. 

Since $X$ has tight hexagons, let $\varphi:D \ra X$ be a tight disc diagram with $\gamma'$ as boundary. The preimage $\varphi^{-1}(K)$ contains at most $n \cdot N$ simplices since $D$ is tight. There exists a connected component $L$ of the combinatorial neighbourhood of the boundary $\partial \varphi^{-1}(K) \subset D$ which contains both the exit edges $e_K(\gamma_{v,w}), e_K(\gamma_{v',w'})$, since the path $\gamma_{a,b} \cup \gamma_{b,b'} \cup  \gamma_{b',a'}$ is disjoint from $K$. This yields a path between those exit edges of length at most $n\cdot N^2 $ since $D$ is tight, which concludes the proof.
\end{proof}

\begin{exm}\label{example_noSAP}
	Here we present an example of an action of a non-hyperbolic group on a hyperbolic graph, satisfying some of the assumptions of Theorem A. Let $X_{(1)}, X_{(2)}$ be the Cayley graphs  with respect to single generators, of $\mathbb{Z}_{(1)},\mathbb{Z}_{(2)}$ being both isomorphic to $\mathbb{Z}$.
	That is $X_{(1)}, X_{(2)}$ are simplicial lines. Let $X$ be the $2$--skeleton of the join of $X_{(1)}$ and $X_{(2)}$. That is, $X$ is the $2$--skeleton of the flag simplicial complex with $1$--skeleton consisting of $1$--skeleta of $X_{(1)}$ and $X_{(2)}$, and edges of the form $\{ v_1,v_2\}$, for $v_i$ being a vertex of $X_{(i)}$. Observe that $X$ has diameter $2$, hence it is hyperbolic.
	The action of $\mathbb{Z}^2=\mathbb{Z}_{(1)}\times \mathbb{Z}_{(2)}$ on $X$ is induced by
	the usual actions of $\mathbb{Z}_{(i)}$ on their Cayley graphs. The action is cocompact, the stabilizers of vertices are isomorphic to $\mathbb{Z}$ and the stabilizers of other simplices are trivial. The action satisfies the conditions (L1), (L2), (G1), (G2) from Theorem A. Note that the complex $X$ has infinite 
	intervals.
	% and does not satisfy the Small Angle Property.    
\end{exm}

%For an example of a complex violating both the small angle property and the tight hexagons property, consider the suspension of the simplicial line with its natural structure of a triangular complex. There exist geodesic bigons between the apices that require an arbitrarily large number of triangles to be filled.

\subsection{Generalization}
\label{s:general}

In this article, we will be considering the collection of \textit{all} geodesics of a given hyperbolic complex. As non locally finite hyperbolic complexes can have a rather wild combinatorial geometry, it is worth mentioning that everything done in this article would work without any change if, instead of the collection of all geodesics of a given complex, we only consider a sub-family $\cP$ of paths satisfying the following conditions:
%
%	A family $\mathcal P$ of paths in $X^{(1)}$ is \emph{good}  if the following conditions are satisfied:
	\begin{enumerate}
		\item \emph{Uniform quasi-geodesics}:
		there exist constants $\lambda \geqslant 1, \varepsilon \geqslant 0$ such that every path in $\mathcal P$ is a $(\lambda, \varepsilon)$--quasi-geodesic,
		\item \emph{Finite intervals}:
		for any two vertices $x, y$  of $X$, the  number of paths in $\mathcal P$ connecting $x$ and $y$ is \emph{finite} and non-zero,
		\item \emph{Stability under restriction}: 
		any subpath of a path in $\mathcal P$ is again in $\mathcal P$.
		\item \emph{Equivariance}:
		the family $\mathcal P$ is $G$--invariant, that is, for every $p\in \mathcal P$ and $g\in G$, the path $g\cdot p$ is again in $\mathcal P$,
		\item \emph{Compatibility with the Gromov boundary}: 
		for every vertex $v$ of $X$ and every point $x$ of the Gromov boundary $\partial X$ of $X$, there exists an infinite path in $\cP$ that is a quasi-geodesic from $v$ to $x$, %Alex: We might need such a path between any two points of $\partial X$, to be checked...
	\end{enumerate}
and by adapting all the definitions used in this article (geodesic hexagons, convexity, the Small Angle Property, etc.) to the setting where only geodesics in $\cP$ are being considered. 

Although our techniques are not meant to apply to that particular case, let us mention that the case of the curve complex of a hyperbolic surface provides such an example of a non locally finite hyperbolic complex with a wild collection of geodesics, but admitting such a nicer sub-family of combinatorial geodesics, namely the collection of \textit{tight geodesics} introduced by Bowditch \cite{BowditchTightGeodesics}. In particular, the first author proved in \cite[Lemma 2.10]{MartinTame} a weakening of the Small Angle Property for tight geodesics. 

%dam vertices v,u,w, simplices \sigma, paths \gamma

%\begin{notation}
%	Paths of $\cP$ are by definition unoriented and we will generally speak of a $\cP$-path \emph{between} two points. If we need a $\cP$-path to be given an orientation, we will speak of an \emph{oriented} $\cP$-path \emph{from} a point \emph{to} another point.
%\end{notation}
%
%
%In what follows, $X$ is a hyperbolic complex together with a good system $\cP$ of combinatorial paths.
%
%\begin{definition}[$\cP$-convexity]
%	A subgraph $K$ of $X^{(1)}$ is said to be $\cP$-convex if every $\cP$-path between two vertices of $K$ is fully contained in $K$.
%\end{definition}

%% file: bdry.tex
\section{Construction of the boundary and the combination theorem}
\label{s:constr}

From now on, we fix an action of a group $G$ on a complex $X$ with quotient space $Y$ satisfying the hypotheses of Theorem A, and let  $G(\cY)= (G_\sigma, \psi_{\sigma, \sigma'}, g_{\sigma, \sigma', \sigma''})$ be a complex of groups over $Y$ associated to that action (see \cite[Section III.$\cC.2.9$]{BridsonHaefliger} for standard results about complexes of groups associated to an action).

\subsection{Construction of the boundary}
\label{s:constrbdry}

In this section, we define the boundary $\partial G$ of the group $G$ and study properties of important subcomplexes associated to points in the boundary of vertex stabilisers. This construction is the analogue of the construction  explained in \cite[Section 2]{MartinBoundaries}, using the language of complexes of groups, and we refer  the reader to \cite[Section 2]{MartinBoundaries} for background on these notions.

% Let $G(\cY)=(G_{\sigma}, \psi_{\sigma, \sigma'}, g_{\sigma, \sigma', \sigma''})$ be a developable complex of groups over a finite polyhedral complex $Y$ with fundamental group $G$.
Recall that the complex of groups $G(\cY)$ consists of the data of \emph{local groups} $G_{\sigma}$ for $\sigma \subset Y$, injective \emph{local maps} $\psi_{\sigma, \sigma'}: G_{\sigma'} \ra G_{\sigma}$ for $\sigma \subset \sigma'$, and  \emph{twist coefficients } $g_{\sigma, \sigma', \sigma''}$ for every $\sigma \subset \sigma' \subset \sigma'' $ subject to the compatibility condition  $  \mbox{Ad}(g_{\sigma, \sigma', \sigma''}) \psi_{\sigma, \sigma''} = \psi_{\sigma, \sigma'}  \psi_{\sigma', \sigma''},$
together with an extra \emph{cocycle condition} that we do not recall here.

Let $F: G(\cY) \ra G$ be a morphism from the complex of groups to $G$ which induces an isomorphism between $G$ and the fundamental group of $G(\cY)$. This amounts to a collection of homomorphisms, still denoted $F$, $G_v \ra G$ together with coefficients $F([\sigma \sigma'])$ for every inclusion $\sigma \subset \sigma'$, satisfying some conditions that we do not recall here. The complex $X$ is isomorphic to the \textit{universal cover} of $G(\cY)$ associated to that morphism, which  is defined as follows: 

$$ \bigg( G \times \coprod_{\sigma \subset Y} \sigma  \bigg) / \simeq$$
where
$$ (gF(g'), x) \simeq (g,x)  \mbox{ if } x \in \sigma, g' \in G_\sigma, g \in G,$$
$$(g,  i_{\sigma, \sigma'}(x) )\sim (g F\big( [\sigma \sigma']\big)^{-1},  x) \mbox{ if } g \in G, \sigma \subset \sigma', x \in \sigma \mbox{ and } i_{\sigma, \sigma'}: \sigma \hra \sigma' \mbox{ is the inclusion.}$$

We recall the following result of \cite[Theorem 2]{MartinCombinationEG}:

\begin{deflem}
 We can associate to $G(\cY)$ the following data:
 \begin{itemize}
  \item for every $\sigma \subset Y$, a (hyperbolic) polyhedral complex endowed with a proper and cocompact action of $G_\sigma$. We denote such a choice of complex by $EG_{\sigma}$, and by $\cV G_{\sigma}$ its vertex set.
  \item for every inclusion $\sigma \subset \sigma'$, a $\psi_{\sigma, \sigma'}$-equivariant polyhedral (quasi-isometric) embedding $\phi_{\sigma, \sigma'}:EG_{\sigma'} \ra EG_{\sigma} $ such that 
  $$\mbox{for every inclusion } \sigma \subset \sigma' \subset \sigma'', \mbox{ we have }    g_{\sigma, \sigma', \sigma''} \circ \phi_{\sigma, \sigma''} = \phi_{\sigma, \sigma'}  \phi_{\sigma', \sigma''}. $$
 \end{itemize}
 We will still denote by $\phi_{\sigma, \sigma'}: \bd G_{\sigma'} \ra \bd G_{\sigma}$ the extension to the Gromov boundaries.\qed
\end{deflem}

\begin{definition}
 We define the space 
$$\cV G =  \bigg( G \times \coprod_{\sigma \subset Y} (\{ \sigma \}\times \cV G_\sigma) \bigg) / \simeq$$
where
$$ (gF(g'), \{\sigma\}, x) \simeq (g,\{ \sigma\}, g'x)  \mbox{ if } x \in \cV G_\sigma, g' \in G_\sigma, g \in G.$$

The canonical projection $G \times \coprod_{\sigma \subset Y} ( \{\sigma\} \times \cV G_\sigma) \ra G \times \coprod_{\sigma \subset Y}  \{\sigma\} $ yields a map 
from $\cV G$ to the vertex set of the first barycentric subdivision of $X$, which we denote $p$. The action of $G$ on $G \times \coprod_{\sigma \subset Y} ( \{\sigma\} \times EG_\sigma)$ on the first factor by left multiplication yields an action of $G$ on $\cV G$, making the projection map $p$ a $G$-equivariant map. 

For every simplex $\sigma$ of $X$, the preimage under $p$ of the barycentre of $\sigma$ is exactly  $\cV G_\sigma $. For an inclusion $\widetilde{\sigma} \subset \widetilde{\sigma}'$ of simplices of $X$ which are lifts of simplices $\sigma \subset \sigma'$ of $Y$, we denote by $\phi_{\widetilde{\sigma}, \widetilde{\sigma}'}: \cV G_{\widetilde{\sigma}'} \ra \cV G_{\widetilde{\sigma}}$ the map sending a point of the form $[g, \{\sigma'\}, x]$ to $[gF\big( [\sigma \sigma']\big)^{-1}, \{\sigma\}, \phi_{\sigma, \sigma'}(x)]$.

% The action of $G$ on $G \times \coprod_{\sigma \in V(\cY)} ( \{\sigma\} \times EG_\sigma)$ on the first factor by left multiplication yields an action of $G$ on $\cV G$, making the projection map $p$ a $G$-equivariant map. 
% 
% Note that $Cl_{D(\cY)}$ can be seen as a complex of spaces over $X$, the fibre of a simplex $[g, \sigma]$ being the classifying space $EG_\sigma$. Indeed, for en edge $[g,a]$ of the first barycentric subdivision of $X$, the gluing map $\phi_{[g\iota_T(a)^{-1}, i(a)], [g, t(a)]}: EG_{i(a)} \ra EG_{t(a)}$ is defined as $\phi_{i(a), t(a)}$.
% 
%  For every simplex $\sigma$ of $X$, we denote by $EG_\sigma$ the fibre over $\sigma$ of that complex of space. For simplices $\sigma, \sigma'$ of $X$ such that $\sigma' \subset \sigma$, we denote by $\phi_{\sigma', \sigma}: EG_\sigma \ra EG_{\sigma'}$ the associated gluing map.
\end{definition} 

\begin{rmk}
The space $\cV G$ is the vertex set, i.e. the $0$-skeleton, of the CW-complex $EG$ constructed in \cite[Section 2]{MartinBoundaries}.
\end{rmk}

\begin{definition}
 We define the space 
$$\Omega G =  \bigg( G \times \coprod_{\sigma \subset Y} ( \{\sigma\} \times \partial G_\sigma) \bigg) / \simeq$$
where

$$ (gF(g'), \{\sigma\}, \xi) \simeq (g, \{\sigma\}, g'\xi)  \mbox{ if } \xi \in \partial G_\sigma, g' \in G_\sigma, g \in G.$$

The canonical projection $G \times \coprod_{\sigma \subset Y} ( \{\sigma\} \times \bd G_\sigma) \ra G \times \coprod_{\sigma \subset Y}  \{\sigma\}$ yields a map 
from $\Omega G$ to the vertex set of the first barycentric subdivision of $X$, which we still denote $p$.

We now define

$$\partial_{Stab} G = \Omega G/ \sim $$
where $\sim$ is the equivalence relation generated by the following identifications:
$$\bigg[g,  \{\sigma\} , \xi\bigg] \sim \bigg[ g F\big( [\sigma \sigma']\big)^{-1},  \{\sigma'\}, \phi_{\sigma,\sigma'}(\xi) \bigg] \mbox{ if } g \in G, \sigma\subset \sigma' \mbox{ and } \xi \in \partial G_\sigma.$$

The action of $G$ on $G \times \coprod_{\sigma \subset Y} ( \left\lbrace \sigma \right\rbrace \times \partial G_\sigma)$ by left multiplication on the first factor yields an action of $G$ on $\Omega G$ and on $\partial_{Stab} G$.

For every simplex $\sigma$ of $X$,  the preimage under $p$ of the barycentre of $\sigma$ is exactly $\bd G_\sigma $. For an inclusion $\widetilde{\sigma} \subset \widetilde{\sigma}'$ of simplices of $X$ which are lifts of simplices $\sigma \subset \sigma'$ of $Y$, we still denote by $\phi_{\widetilde{\sigma}, \widetilde{\sigma}'}: \bd G_{\widetilde{\sigma}'} \ra \bd G_{\widetilde{\sigma}}$ the map extended to the Gromov boundaries, sending a point of the form $[g, \{\sigma'\}, \xi]$ to $[gF\big( [\sigma \sigma']\big)^{-1}, \{\sigma\}, \phi_{\sigma, \sigma'}(\xi)]$.
\end{definition} 

We recall the following result of \cite[Proposition 4.4]{MartinBoundaries}. Note that the proof therein still holds in our case, as it does not use the geometry of the complex but only the fact that loops in fixed-point sets are contractible by condition (G2). We point out that this is the only result of this article where condition (G2) is used.

\begin{deflem}
 Let $\sigma$ be a simplex of $X$. Then the projection $\pi : \Omega G \ra \bds$ is injective on $\partial G_{\sigma}$. We thus still denote by $\bd G_{\sigma}$ the image of $\bd G_{\sigma} \subset \Omega G$ in $\bds$. \qed
\label{xiloop}
\end{deflem}

% \begin{proof}
% This is essentially Proposition 4.4 of \cite{MartinBoundaries}. The only argument in the proof which is not purely combinatorial and must be modified is the construction of the so-called \textit{hull} of a loop contained in a domain \cite[p.???]{MartinBoundaries}. This was proved by using the fact that the pointwise stabiliser of such a loop has a simply-connected fixed-point set, which is condition ??? in our definition.
% \end{proof}

\begin{definition}
 We define the \textit{boundary} $$\partial G := \partial_{Stab} G \sqcup \partial X$$ and the \textit{compactification} $$\cC G := \cV G \sqcup \partial G.$$
 For every simplex $\sigma$ of $X$, we define the subset 
 $$\cC G_{\sigma} := \cV G_{\sigma} \sqcup \partial G_{\sigma}.$$
\end{definition}

\begin{rmk}\label{rmk:explanation}
In 	\cite{MartinBoundaries}, a slightly different compactification of $G$ was considered. Indeed, instead of looking at $\cV G \sqcup \partial G,$ the first author considered the space $EG \sqcup \partial G,$ where $EG$ was a classifying space for proper actions of $G$. The reason for this was to show that the resulting boundary yielded an $E\cZ$-structure in the sense of Farrell--Lafont, a structure that implies the Novikov conjecture for the group. As in this article we are only interested in hyperbolic groups, we prefer to introduce a very close compactification, but which has the advantage of being easier to handle.
\end{rmk}

In the next section, we will define a topology on  $\cC G$ that will turn it into a compact space, justifying the terminology.

\subsection{Definition of the topology}
\label{s:topology}

We now define a topology on the compactification $\cC G$ following \cite{MartinBoundaries}, by defining a basis of neighbourhoods at every point. Note that in \cite{MartinBoundaries}, a basis of \textit{open} neighbourhoods was defined. In this article, considering the slightly different compactification used for the reasons outlined in Remark \ref{rmk:explanation},  it will be easier in our case to consider (not necessarily open) neighbourhoods, which explains the slight changes between the definitions.

Although we are mostly interested in the topology of the boundary $\partial G$, it will be important to have a topology on $\CG$ in some of the proofs. Since points of $\cC G$ are of three different kinds ($\cV G$, $\partial X$ and $\partial_{Stab} G$), we treat these cases separately.

Recall that by hypothesis, the system of geodesic paths in the universal cover $X$ of $G(\cY)$ satisfies the Small Angle Property.

\begin{definition}[based geodesic]
We fix once and for all a base-vertex $v_0$ of $X$. An oriented geodesic $\gamma_{v,w}$ between two generalised vertices of $\overline{X}$ is said to be \emph{based} if $v=v_0$. In such a case, we will simply denote it $\gamma_w$. 
\end{definition}

\subsubsection{Domains}

\begin{definition}(domains and projection)
We define \emph{domains} of points of $\cC G$ as follows:
\begin{itemize} \item Let $\xi \in \partial_{Stab} G$. We define the \textit{domain} of $\xi$, denoted $D(\xi)$, as the  subgraph of the $1$-skeleton of $X$ spanned by vertices $v$ such that $\xi \in \partial G_v$. 
%We denote by $V(\xi)$ the set of vertices of $D(\xi)$.
 \item Let $\eta \in \partial X$. We define the domain of $\eta$ as the singleton $\{\eta\} \subset \overline{X}$.
 \item Let $v$ be a point of $\VG$. We define the domain of $v$ as the unique simplex $\sigma$ of $X$ such that $v \in \cV G_\sigma$.
% \item Let $v$ be a point of $\VG$ and $\sigma$ the unique simplex of $X$ such that $v \in \cV G_\sigma$. We define the domain of $v$ as the singleton $\{\sigma\} \subset X'$, seen as a vertex of the first barycentric subdivision of $X$.
 \end{itemize}
 We extend the projection map $p$ from $\cV G$ to the vertex set of the first barycentric subdivision of $X$ to a (coarse) projection map from $\cC G$ to the set of subspaces of $\overline{X}$ by sending every point of $\cC G$ to its domain.
\label{def:domain}
\end{definition}

The following proposition is crucial in defining a topology on $\cC G$. It is also the unique moment in this article where the weak acylindricity condition (G1) is used.

\begin{prop}
Domains are combinatorially convex subcomplexes of $X$ with a uniformly bounded above number of edges.
\label{finitedomain}
\end{prop}

\begin{proof}
This is essentially Proposition 4.2 of \cite{MartinBoundaries} (in the original statement, the uniform bound is not mentioned, though it follows immediately from the proof). Note in particular that the fact that domains are bounded is a direct consequence of the weak acylindricity (G1). The only places where the arguments must be adapted are the following:

\textit{Combinatorial convexity}: The proof of combinatorial covexity of domains used the uniqueness of CAT(0) geodesics. By contrast, here we only known that intervals are finite. However, the proof does extend to this setting. Indeed, for vertices $v,v'$ in a domain, let $H$ be the pointwise stabiliser of the pair $\{v,v'\}$. By the finite interval condition, there are only finitely many geodesics between $v$ and $v'$, so there exists a finite index subgroup $H'$ of $H$ which pointwise stabilises every geodesic between $v$ and $v'$. As taking finite index subgroups does not change limit sets, the proof of \cite[Lemma 4.7]{MartinBoundaries} works with  $H'$ instead of $H$ itself, and thus every geodesic between $v$ and $v'$ is  contained in the domain $D(\xi)$ .

\textit{Finiteness of domains}: The proof that domains are locally finite used the so-called \textit{Limit Set property} and \textit{Finite Height property} (see \cite[Proposition 4.2]{MartinBoundaries}). The fact that such properties are satified for complexes of groups with hyperbolic local groups and quasi-convex embeddings as local maps was proved in \cite[Lemmas 9.4 and 9.7]{MartinBoundaries}.
\end{proof}

\begin{definition}\label{def:dmax}
 We denote by $d_{max}$ a uniform bound on the number of edges of domains of points of $\bds$.
\end{definition}

\begin{notation}
Let $\gamma$ be a based geodesic, $\xi$ a point of $\partial_{Stab} G$,  and suppose that $\gamma$ goes through $D(\xi)$. We will simply denote $e_\xi(\gamma)$ the exit edge associated to the domain $D(\xi)$.
\end{notation}

\subsubsection{Neighbourhoods of a point of $\cV G$}

\begin{definition}
 Let $v \in \cV G$. We define a collection of neighbourhoods $\cB_{\CG}(v)$ of $v$ in $\cC G$ as the family of finite subsets of $\cV G$ containing $v$.
\end{definition}

\subsubsection{Neighbourhoods of a point of $\partial X$}

We now turn to the case of points of the boundary of $X$. Recall that since $X$ is a hyperbolic complex with countably many simplices, the bordification $\overline{X} = X \cup \partial X$ obtained by adding the Gromov boundary of $X$ has a natural metrisable topology, though not necessarily compact if $X$ is not locally finite. For every $\eta \in \partial X$, a basis of neighbourhoods of $\eta$ in that bordification is given by the family of
$$W_k(\eta) = \left\{x \in \overline{X} \mbox{ such that } \langle x, \eta \rangle_{v_0} \geqslant k \right\} ,$$

where $\langle \cdot, \cdot \rangle_{v_0}$ denotes the Gromov product with respect to $v_0$.
We denote by $\cB_{\overline{X}}(\eta)$ this basis of neighbourhoods of $\eta$ in $\overline{X}$. Endowed with that topology, $\overline{X}$ is a second countable metrisable space.
\iffalse
Note that the topology of $\overline{X}$ satisfies the following property:

\begin{lem}
Let $\eta \in \partial X$, $U$ a neighbourhood of $\eta$ in $\overline{X}$ and $k \geqslant 0$. Then there exists a neighbourhood $U'$ of $\eta$ in $\overline{X}$ that is contained in $U$ and such that $d(U' \cap X, X \setminus U) >k.$
\label{hyperbolicnesting}
\end{lem}

\begin{proof}
The definition of the topology of $\overline{X}$ implies the following: if $(x_n)$ and $(y_n)$ are two sequences of $X$ such that $d(x_n,y_n)$ is bounded, then $(x_n)$ converges to a point of $\partial X$ if and only if $(y_n)$ converges to the same point. Reasoning by contradiction thus implies the lemma.
\end{proof}
\fi 

\begin{definition}
 Let $\eta \in \partial X$, and let $U$ be a neighbourhood of $\eta$ in $\cB_{\overline{X}}(\eta)$. We define a neighbourhood $V_U(\eta)$ as the set of elements $z \in \cC G$ whose domain is contained in $U$. 
When $U$ runs over the basis $\cB_{\overline{X}}(\eta)$ of neighbourhoods of $\eta$ in $\overline{X}$, this defines a collection of neighbourhoods of $\eta$, which we denote $\cB_{\cC G}(\eta)$.
\end{definition}

\subsubsection{Neighbourhoods of a point of $\partial_{Stab} G$}

We finally define neighbourhoods for points in $\partial_{Stab} G$. Since, in $\partial_{Stab} G$, boundaries of stabilisers of vertices are glued together along boundaries of stabilisers of edges, we will construct neighbourhoods in $\cC G$ of a point $\xi \in \partial_{Stab} G$ using neighbourhoods of the representatives of $\xi$ in the various $\cC G_v$, where $v$ runs over the vertices of the domain of $\xi$.

% \begin{definition}[the convergence property]
%  We say that an $E\cZ$-complex of spaces compatible with $G(\cY)$ satisfies the \textit{convergence property} if, for every pair of simplices $\sigma \subset \sigma'$ in $Y$ and every injective sequence $(g_nG_{\sigma'})$ of cosets of $G_\sigma / G_{\sigma'}$, there exists a subsequence such that $(g_{\varphi(n)}\overline{EG_{\sigma'}})$ uniformly converges to a point in $\overline{EG_\sigma}$.
% \label{convergenceproperty}
% \end{definition}

\begin{definition}[$\xi$-family]
 Let $\xi$ be a point of $\partial_{Stab} G$. A collection $\cU$ of open sets $\left\lbrace U_\sigma, \sigma \subset D(\xi) \right\rbrace $ is called a $\xi$-\textit{family} if for every pair of vertices $v,v'$ of $D(\xi)$ joined by an edge $e$, and for every $x \in \cC G_e$, we have
$$\phi_{v,e}(x) \in U_v \Leftrightarrow \phi_{v',e}(x) \in U_{v'}.$$
\label{defxifamily}
\end{definition}

It was proved in  \cite[Proposition 4.12]{MartinBoundaries} that for every point $\xi \in \partial_{Stab} G$ and every collection $(U_v)_{v \in D(\xi)}$, where each $U_v$ is a neighbourhood of $\xi$ in $ \VG_\sigma$, there exists a $\xi$-family $\cU'$ such that $U_v' \subset U_v$ for every vertex $v$ of $D(\xi)$. As the proof did not use any CAT(0) geometry, the same holds in our situation.
%\label{xifamilies}
%\end{lem}

\begin{definition}[Cone] Let $\xi \in \partial_{Stab} G$ and $\cU$ be a $\xi$-family.
We define the \textit{cone} $\mbox{Cone}_{\cU}(\xi)$ as the set of generalised vertices $v$ of $X \cup \partial X$ such that every based geodesic $\gamma_v$ goes through $D(\xi)$ (in particular, $v \notin D(\xi)$) and such that for the unique vertex $w$ of $D(\xi) \cap e_{\xi}(\gamma_v)$, we have the following inclusion in $\cC G_w$: 
$$\cC G_{e_{\xi }(\gamma_v)} \subset U_w.$$
\end{definition}

\begin{definition}
Let $\xi \in \partial_{Stab} G$ and $\cU$ be a $\xi$-family. We define the neighbourhood $V_{\cU}(\xi)$ of $\xi$ as the set of points $z \in \CG$ such that $D(z) \setminus D(\xi) \subset \mbox{Cone}_{\cU}(\xi)$ and such that for every vertex $v$ of $D(z) \cap D(\xi)$ we have $z \in U_v$.

This collection of neighbourhoods of $\xi$ in $\CG$, for $\cU$ ranging over all possible $\xi$-families, is denoted $\cB_{\cC G}(\xi)$. 
\label{definitionneighbourhoodsxi}
\end{definition}

\begin{rmk}
Since this definition involves \emph{based} geodesics, these neighbourhoods depend on the chosen basepoint $v_0$. In the Appendix, we will need to allow change of basepoints, in which case we will indicate the basepoint used to define the various cones and neighbourhoods  $\mbox{Cone}_{\cU}(\xi)$, $V_{\cU}(\xi)$ as a superscript. In that case, we will speak of the topology (of $\CG$) \textit{centred} at a given vertex.
\end{rmk}

\begin{definition}
 We define a topology on $\CG$ by taking the topology generated by the elements of $\cB_{\CG}(x)$, for every $x \in \CG$. That is, a subset $U$ of $\cC G$ is open if for every point $x \in U$, there exists a neighbourhood $U' \in \cB_{\cC G}(x)$ contained in $U$.
 We denote by $\cB_{\CG}$ the set of elements of $\cB_{\CG}(x)$ when $x$ runs over $\CG$. 
\end{definition}

\subsection{Uniform convergence group actions and the combination theorem}
\label{s:conv}

 To prove the Combination Theorem A, we use a topological characterisation of hyperbolicity due to Bowditch \cite{BowditchTopologicalCharacterization}, using the boundary $\partial G$, together with the topology defined in the previous section, as a candidate for the Gromov boundary of $G$.

Recall that a group $\Gamma$ acting on a compact metrisable space $M$ with more than two points is  a \textit{convergence group} if the following holds: For every injective sequence $(\gamma_n)$ of elements of the group, there exist, up to taking a subsequence of $(\gamma_n)$, two points $\xi_+, \xi_-$ of $M$  such that for every compact subspace $K \subset M \setminus \left\lbrace \xi_- \right\rbrace$, the sequence of translates $(\gamma_{n}K)$  uniformly converges to $\xi_+$.

A hyperbolic group $\Gamma$ is always a convergence group on $\Gamma \cup \partial \Gamma$ (see for instance \cite{FredenConvergenceGroup}). %As noted in \cite[???]{MartinBoundaries}, it follows immediately that 
It also is a convergence group on $\overline{E\Gamma}$.

Recall that for a group $\Gamma$ that is a convergence group on a compact metrisable space $M$ with at least three points, a point $\zeta$ in $M$ is called a \textit{conical limit point} if there exists a sequence $(\gamma_n)$ of elements of $\Gamma$ and two points $\xi_- \neq \xi_+$ in $M$, such that $\gamma_n \zeta \ra \xi_-$ and $\gamma_n \zeta' \ra \xi_+$ for every $\zeta' \neq \zeta$ in $M$. The group $\Gamma$ is a \textit{uniform convergence group} on $M$ if every point of $M$ is a conical limit point. By a celebrated result of Bowditch \cite{BowditchTopologicalCharacterization}, a group $\Gamma$ that is a uniform convergence group on a compact metrisable space $M$ with more than two points is hyperbolic, and $M$ is $\Gamma$-equivariantly homeomorphic to the Gromov boundary  of $\Gamma$. 

Our main result is the following theorem implying immediately Theorem A from Introduction.

\begin{thm}\label{thm:combination}
 The boundary $\partial G$ is a compact metrisable space, on which the group $G$ acts as a uniform convergence group. In particular, $G$ is hyperbolic and $\partial G$ is $G$-equivariantly homeomorphic to the Gromov boundary of $G$. Moreover, every local group of $G(\cY)$ embeds in $G$ as a quasi-convex subgroup. \qed
\end{thm}

Theorem \ref{thm:combination} extends the results of \cite{MartinBoundaries} to the case of a group acting on a polyhedral complex endowed with a geometry that is non-positively curved in a combinatorial sense. The proofs, which rely on properties reminiscent of non-positive curvature, are very close to the original proofs of \cite[Section 9]{MartinBoundaries} in the CAT(0) setting, and are given in the Appendix.

%% file: wsys.tex
\section{Examples of group actions on combinatorially nonpositively curved complexes}
\label{s:examples}

\subsection{Combinatorial geometry of  weakly systolic complexes}
\label{s:wsys}

The main results of this subsection -- Proposition~\ref{p:sap-wsys}, Proposition~\ref{p:fixpt-wsys}, and Proposition  \ref{p:finint-wsys} -- allow to deduce
Theorem B in Introduction, from Theorem A. 
Throughout this subsection we do not assume the subcomplexes to be hyperbolic. 

\begin{prop}
\label{p:sap-wsys}
A weakly systolic complex has tight hexagons. In particular, it satisfies the Small Angle Property.
\end{prop}
\begin{prop}
\label{p:fixpt-wsys}
For any action on a weakly systolic complex loops in fixed-point sets are contractible.
\end{prop}

Recall that a simplicial complex is \emph{flag} if every set of pairwise adjacent vertices is a simplex.
A flag simplicial complex \emph{without infinite simplices} is a complex without infinite cliques 
(complete graphs) in its $1$--skeleton.

\begin{prop}
\label{p:finint-wsys}
Weakly systolic complexes without infinite simplices have finite intervals.
\end{prop}

\begin{definition}[weakly systolic complex]
\label{d:wsys}
A flag simplicial complex $X$ is \emph{weakly systolic} if for every vertex $v\in X$ and for every $n\geqslant 1$
the following two conditions are satisfied: 
\begin{description}
  \item[({\bf E}) \rm(edge condition)] For every
edge $e$ with both end-vertices at distance $n$ from $v$, there is a vertex $u$ at distance $n-1$ from $v$, adjacent to the end-vertices of $e$.
  \item[({\bf V}) \rm{(vertex condition)}] For every vertex $w$ at distance $n$ from $v$, the set of vertices adjacent
to $w$ and at distance $n-1$ from $v$ induces a clique (full graph), that is, they are all adjacent.
\end{description}
\end{definition}
\medskip

\noindent
{\bf Examples.} An important class of examples of weakly systolic complexes are \emph{systolic complexes}.
These are simply connected flag simplicial complexes in which every loop (that is, a closed path) of length
(number of edges) $4$ or $5$ has a \emph{diagonal}, that is, an edge connecting non-consecutive vertices.
Simplicial trees are systolic. A basic example of an infinite systolic complex is the triangulation of
the Euclidean plane by equilateral triangles. Other important examples are highly dimensional hyperbolic 
pseudomanifolds acted geometrically upon groups constructed in \cite{SystolicComplexes}.
Systolic complexes have a particular large-scale geometry -- they `do not contain asymptotically' spheres of dimension two and more; see e.g.\ \cite{SystolicComplexes,O-cnpc,ChO}. An example of a weakly systolic complex of large-scale geometry different than the one of systolic complexes is obtained as follows. Let $X$ be a CAT($-1$) cubulation of the real $n$--dimensional
Lobachevski hyperbolic space $\mathbb H^n$, for $n=3,4$. Let $Th(X)$ be the \emph{thickening} of $X$ (see e.g.\ \cite{O-cnpc,ChO}), that is, a flag simplicial complex with the set of vertices being $X^{(0)}$ and two vertices connected by an edge iff they are contained in a common cube of $X$. Then $Th(X)$ is weakly systolic but not quasi-isometric to a systolic complex. More generally, the thickening of every CAT($-1$) cubical complex is weakly systolic. For further examples of weakly systolic complexes and their automorphism groups see e.g.\ \cite{SystolicComplexes,O-cnpc,ChO}.
\medskip

%The \emph{interval}
%$I(u,v)$ between $u$ and $v$ consists of all vertices on shortest
%$(u,v)$-paths, that is, of all vertices (metrically) \emph{between} $u$
%and $v$: $I(u,v)=\{ x\in V: d(u,x)+d(x,v)=d(u,v)\}$.

Weak systolicity can be seen as a combinatorial analogue of a metric non-positive curvature.
In particular, in \cite{O-cnpc,ChO} a local-to-global characterization of weakly systolic complexes
is proved: a universal cover of a flag simplicial complex satisfying a \emph{local edge condition} and 
a \emph{local vertex condition} is weakly systolic. `Local' versions of the conditions (E) and (V) above 
are obtained by restricting the values of $n$ to $1,2,3$.
\medskip

We now proceed to the proofs of Propositions~\ref{p:sap-wsys}, \ref{p:fixpt-wsys} and \ref{p:finint-wsys}.
In the remaining part of this subsection we assume that $X$ is a weakly systolic complex.

\begin{lem}(simple bigon filling)
\label{l:simpbigon}
Let $\gamma_{uv}=(z^0:=v,z^1,\ldots,z^n:=u)$ and $\gamma'_{uv}=(z'^0:=v,z'^1,\ldots,z'^n:=u)$ be two geodesics between vertices $v$ and $u$ such that $z^i\neq z'^i$, for $i\neq 0,n$. Then there exists a disc diagram $D \to X$ for the loop $\gamma_{uv}\circ {\gamma'_{uv}}^{-1}$ which is an embedding, with internal vertices degree $6$ and the boundary vertices degree at most $5$.
\end{lem}
\begin{proof}
By Claim 1 in the proof of Theorem~3.1 in \cite{ChO}, there is a minimal disc diagram $D \to X$ with systolic $D$ (i.e.\ every internal vertex has degree at least $6$).
We will denote vertices in $D$ with tilde, like $\wt v$, and their images in $X$ without tilde, that is, $D\to X\colon \wt v\mapsto v$.
In particular,  the paths $\wt{\gamma}_{uv}=(\wt z^0:=\wt v,\wt z^1,\ldots,\wt z^n:=\wt u)$ and $\wt{\gamma}'_{uv}=(\wt z'^0:=\wt v,\wt z'^1,\ldots,\wt z'^n:=\wt u)$ are combinatorial geodesics in $D$ mapped isometrically onto, respectively, $\gamma_{uv}$ and $\gamma'_{uv}$.
The combinatorial Gauss-Bonnet formula reads:
\begin{align}
\sum_{\wt v\in D^{(0)}} {\rm def}(\wt v) = 6,
\end{align}
where the \emph{defect} ${\rm def}(\wt v)$ of a vertex $\wt v$ is equal to $3$ (respectively, $6$) minus the number of triangles containing $\wt v$, for 
$\wt v$ lying on the boundary (respectively, in the interior) of $D$.
From the vertex condition (V) it follows that ${\rm def}(\wt u)={\rm def}(\wt v)=2$. The sum of defects on a geodesic is at most $1$ (see e.g.\ \cite[Fact 3.1]{Elsner2009-flats}), that is
\begin{align}
\begin{split}
{\rm def}(\wt z^1)+{\rm def}(\wt z^2)+\cdots+{\rm def}(\wt z^{n-1})\leqslant 1,\\
%\end{align}
%and
%\begin{align}
{\rm def}(\wt z'^1)+{\rm def}(\wt z'^2)+\cdots+{\rm def}(\wt z'^{n-1})\leqslant 1.
\end{split}
\end{align}
Since all interior vertices have nonpositive defect, we get that every interior vertex has defect $0$ and the sum
of defects on each boundary geodesic is $1$. Thus $D$ is a flat disc (see e.g.\ \cite[Lemma 3.5]{Elsner2009-flats}), that is, a subcomplex of an equilateral triangulation
of the Euclidean plane homeomorphic to a disc.

Therefore, the disc $\triangle$ consists of vertices $\wt z^i_j$, for $i=0,1,\ldots,n$ and $j=0,1,\ldots,k(i)$, satisfying the following conditions (see Figure~\ref{f:bigdiag}):
\begin{enumerate}[(a)]
\item
\label{e:a}
$d(\wt z^i_j,\wt v)=i=n-d(\wt z^i_j,\wt u)$;
\item
$\wt z^i_0=\wt z^i$ and $\wt z^i_{k(i)}=\wt z'^i$;
\item
if $\wt z^i_j$ and $\wt z^{i'}_{j'}$ are adjacent with $i \leqslant i'$ then $i'-i \leqslant 1$ and:
\begin{enumerate}[(i)]
\item
if $i=i'$ then $|j-j'|=1$,
\item
if $i'=i+1$ then $|j-j'|\leqslant 1$;
\end{enumerate}
\item
if $\wt z^i_j$ is adjacent to both $\wt z^{i+1}_{j'}$ and $\wt z^{i+1}_{j''}$ then $|j'-j''|\leqslant 1$.
\end{enumerate}

\begin{figure}[h]
\begin{center}
\includegraphics[scale=0.75]{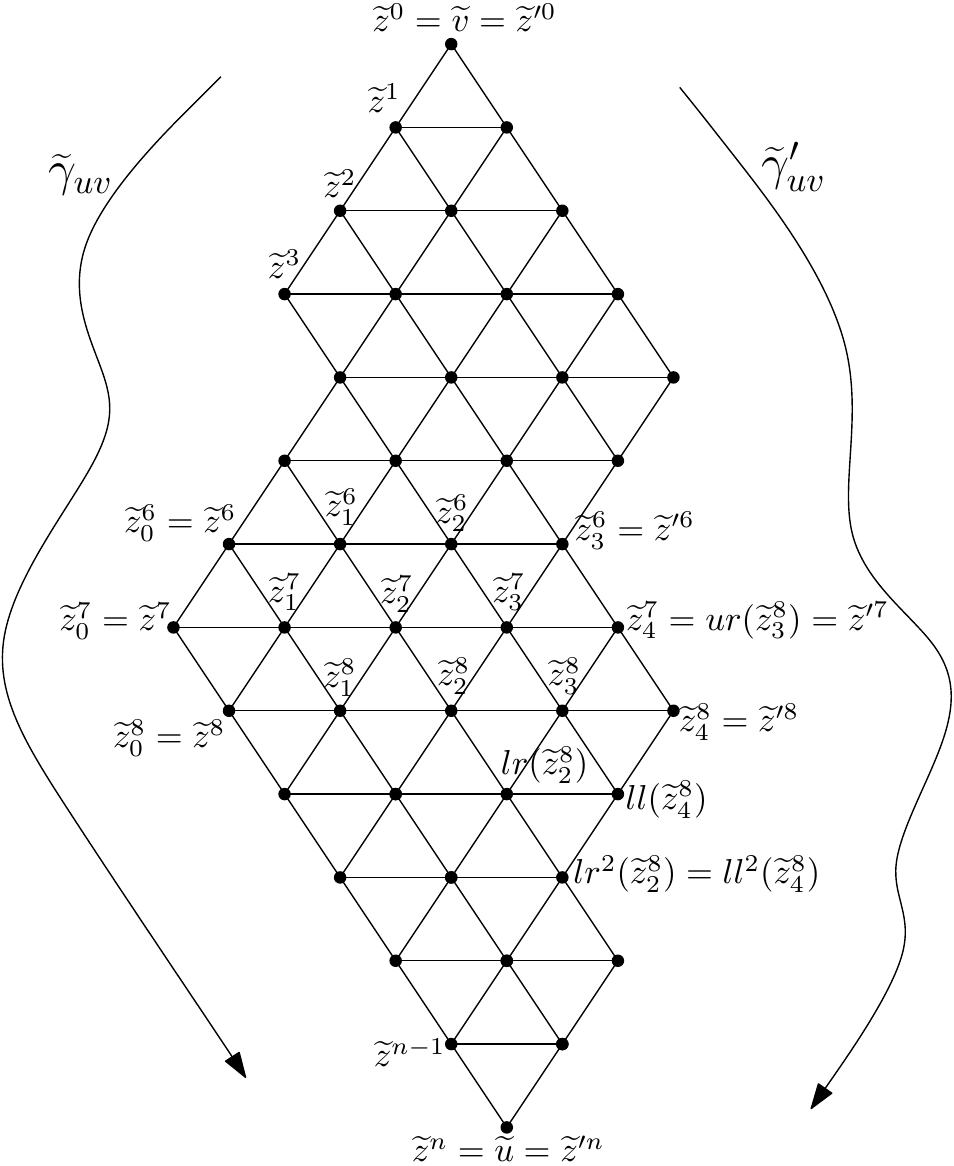}
\end{center}
\caption{Geodesic bigon.}\label{f:bigdiag}
\end{figure}
Observe that the minimal disc diagram for a bigon of length $n$ has area at most $n^2/2$.
Let $ul(\wt z^i_j)$ and $ll(\wt z^i_j)$ (respectively, $ur(\wt z^i_j)$ and $lr(\wt z^i_j)$) denote the vertices $\wt z^{i+1}_k$ and $\wt z^{i-1}_m$ adjacent to $\wt z^i_j$ and with minimal $k$ and $m$ (respectively, maximal $k$ and $m$). (Here $ul$ and $ll$ come from `upper-left' and `lower-left', and so on -- see Figure~\ref{f:bigdiag}). Further, we set $ul^0(\wt z^i_j):=\wt z^i_j$ and, by induction,
$ul^{k+1}(\wt z^i_j):=ul(ul^{k}(\wt z^i_j))$. Similarly we define $ll^k(\wt z^i_j), ur^k(\wt z^i_j)$, and $lr^k(\wt z^i_j)$.

We show now that the disc diagram $D \to X$ is an embedding.
By (\ref{e:a}) we could have that $z^i_j=z^k_l$ only if $i=k$.
Suppose that $z^i_j=z^i_l$, for $j<l$. Then the diagram $D \to X$ can be modified in the following way.
Consider an equilateral triangle $T_l$ in $D$ bounded by the paths: $(\wt z^i_j,\wt z^i_{j+1},\ldots,\wt z^i_{l})$, $\wt \gamma_{lr}:=(\wt z^i_j, lr(\wt z^i_j),lr^{2}(\wt z^i_j)),\ldots, lr^{l-j}(\wt z^i_j))$, 
$\wt \gamma_{ll}:=(\wt z^i_l, ll(\wt z^i_l),ll^{2}(\wt z^i_l)),\ldots, ll^{l-j}(\wt z^i_l))$.
Analogously, we define paths $\wt \gamma_{ur}$, $\wt \gamma_{ul}$, and the triangle $T_u$ using $ur$ and $ul$.
Since $z^i_j=z^i_l$ we can modify $D \to X$ by `squeezing' $(\wt z^i_j,\wt z^i_{j+1},\ldots,\wt z^i_{l})$ to a point $\wt z^i_j$ and
filling the bigon $\gamma_{lr} \gamma_{ll}^{-1}$ by a disc $D_l$, similarly with the bigon $\gamma_{ur} \gamma_{ul}^{-1}$ filled by $D_u$ (where by $\gamma_{ab}^{\pm 1}$ we denote the image of $\wt \gamma_{ab}^{\pm 1}$ via $D\to X$). 
Replacing $T_l, T_u$ by, respectively, $D_l, D_u$  we obtain a new disc diagram $D' \to X$ for the loop $\gamma_{uv} {\gamma'_{uv}}^{-1}$ (see Figure~\ref{f:squeezing}).
\begin{figure}[h]
\begin{center}
\includegraphics[scale=0.75]{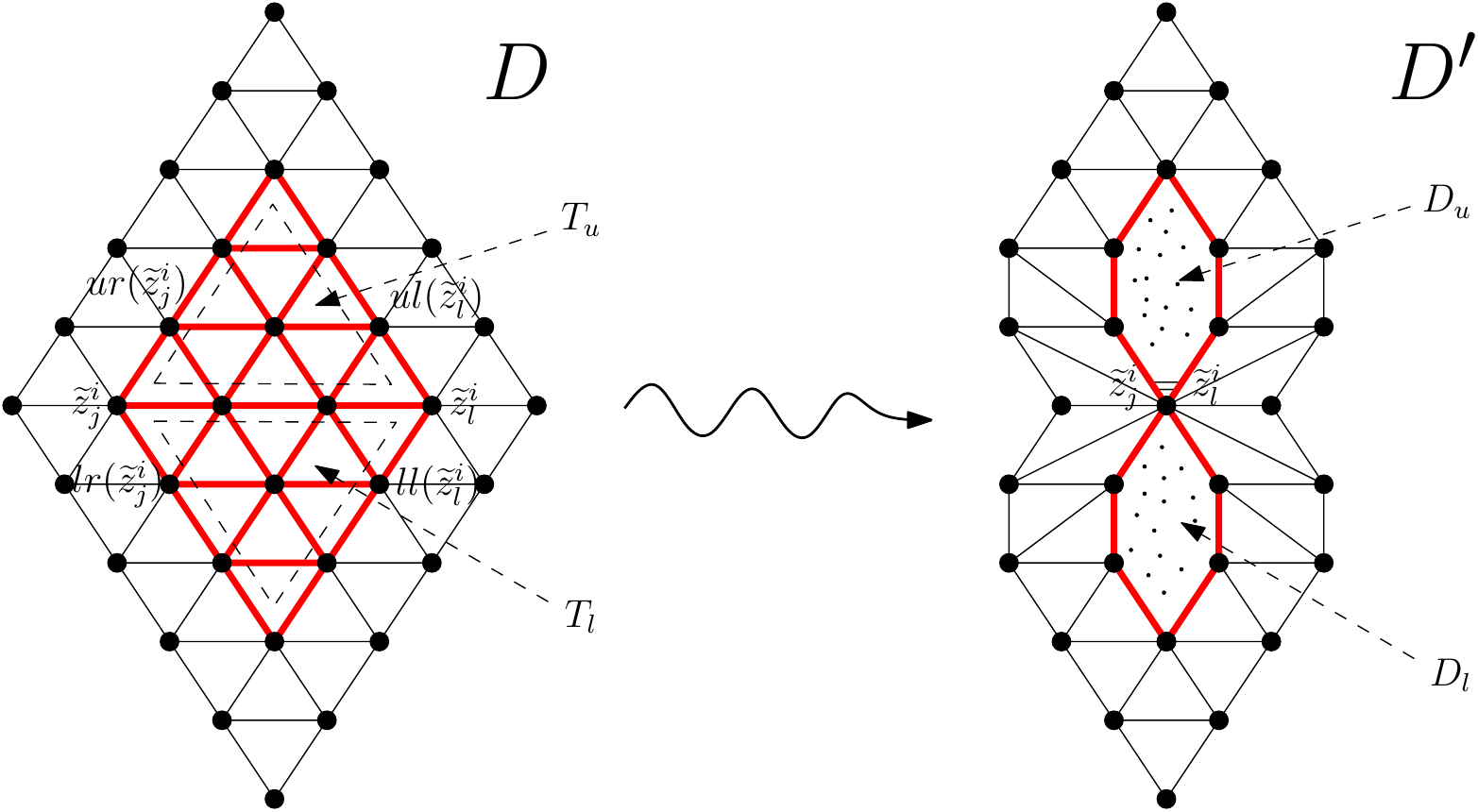}
\end{center}
\caption{From $D \to X$ to $D' \to X$.}\label{f:squeezing}
\end{figure}
Since the area of $T_l$ (and of $T_u$) is $(l-j)^2$, and the area of a minimal disc diagram for
the bigon $\gamma_{lr} \gamma_{ll}^{-1}$ is at most $(l-j)^2/2$ (the same for $\gamma_{ur} \gamma_{ul}^{-1}$), it follows
that the area of $D'$ is strictly less than the area of $D$. This contradicts the minimality of $D$.
Hence $D\to X$ is an embedding.
\end{proof}

\begin{lem}(bigon filling)
\label{l:bigon}
Let $\gamma_{uv}=(z^0:=v,z^1,\ldots,z^n:=u)$ and $\gamma'_{uv}=(z'^0:=v,z'^1,\ldots,z'^n:=u)$ be two geodesics between vertices $v$ and $u$. Then there exists a disc diagram $D \to X$ for the loop $\gamma_{uv}{\gamma'_{uv}}^{-1}$ which is an embedding, with internal vertices degree $6$ and the boundary vertices degree at most $5$.
\end{lem}
\begin{proof}
Decompose $\gamma$ and $\gamma'$ into subgeodesics so that the corresponding loops are simple
and then use Lemma~\ref{l:simpbigon}.
\end{proof}

Recall that one-skeleta of weakly systolic complexes are \emph{weakly modular} graphs \cite{ChO}. 
Three vertices $v_1,v_2,v_3$ of a graph form a {\it metric triangle}
$v_1v_2v_3$ if the intervals $I(v_1,v_2)$, $I(v_2,v_3),$ and
$I(v_3,v_1)$ pairwise intersect only in the common end-vertices. If
$d(v_1,v_2)=d(v_2,v_3)=d(v_3,v_1)=k,$ then this metric triangle is
called {\it equilateral} of {\it size} $k$.

\begin{lem} \cite{Ch_metric} \label{l:metrtr} A graph
  is weakly modular if and only if for any metric triangle $u'v'w'$
   and any two vertices $x,y\in I(v',w'),$ the equality
  $d(u',x)=d(u',y)$ holds. In particular, all metric triangles of weakly modular graphs are equilateral.
\end{lem}

Lemma~\ref{l:metrtr} implies immediately the following.

\begin{lem}(metric triangle filling)
\label{l:metrsys}
Let $u'v'w'$ be a metric triangle in $X$. Then there exists a systolic equilateral triangle $D$
and a disc diagram $D \to X$ which is an isometric embedding and maps the three vertices of $D$ onto $u',v',w'$.
\end{lem}

\begin{prop}(tight geodesic triangle)
\label{p:tightsys}
Let $uvw$ be vertices of a geodesic triangle with sides $\gamma_{uv},\gamma_{vw},\gamma_{wu}$. Then there exists a disc diagram
$D \to X$ for the loop $\gamma=\gamma_{uv}\circ \gamma_{vw} \circ \gamma_{wu}$ which is at most $4$--to--$1$ and
whose every vertex has degree at most $14$. 
\end{prop}
\begin{proof}
Let $u'$ be a vertex in $I(u,w)\cap I(u,v)$ that is at a maximal distance from $u$. Analogously we define $v'$ and $w'$.
Then $u'v'w'$ is a metric triangle -- see Figure~\ref{f:metrtr}. Choose a geodesic $\gamma_{u'v'}$ between vertices $u'$ and $v'$, and similarly choose
geodesics $\gamma_{v'w'}$,$\gamma_{w'u'}$, $\gamma_{uu'}$, $\gamma_{vv'}$, $\gamma_{ww'}$.
\begin{figure}[h]
\begin{center}
\includegraphics[scale=0.65]{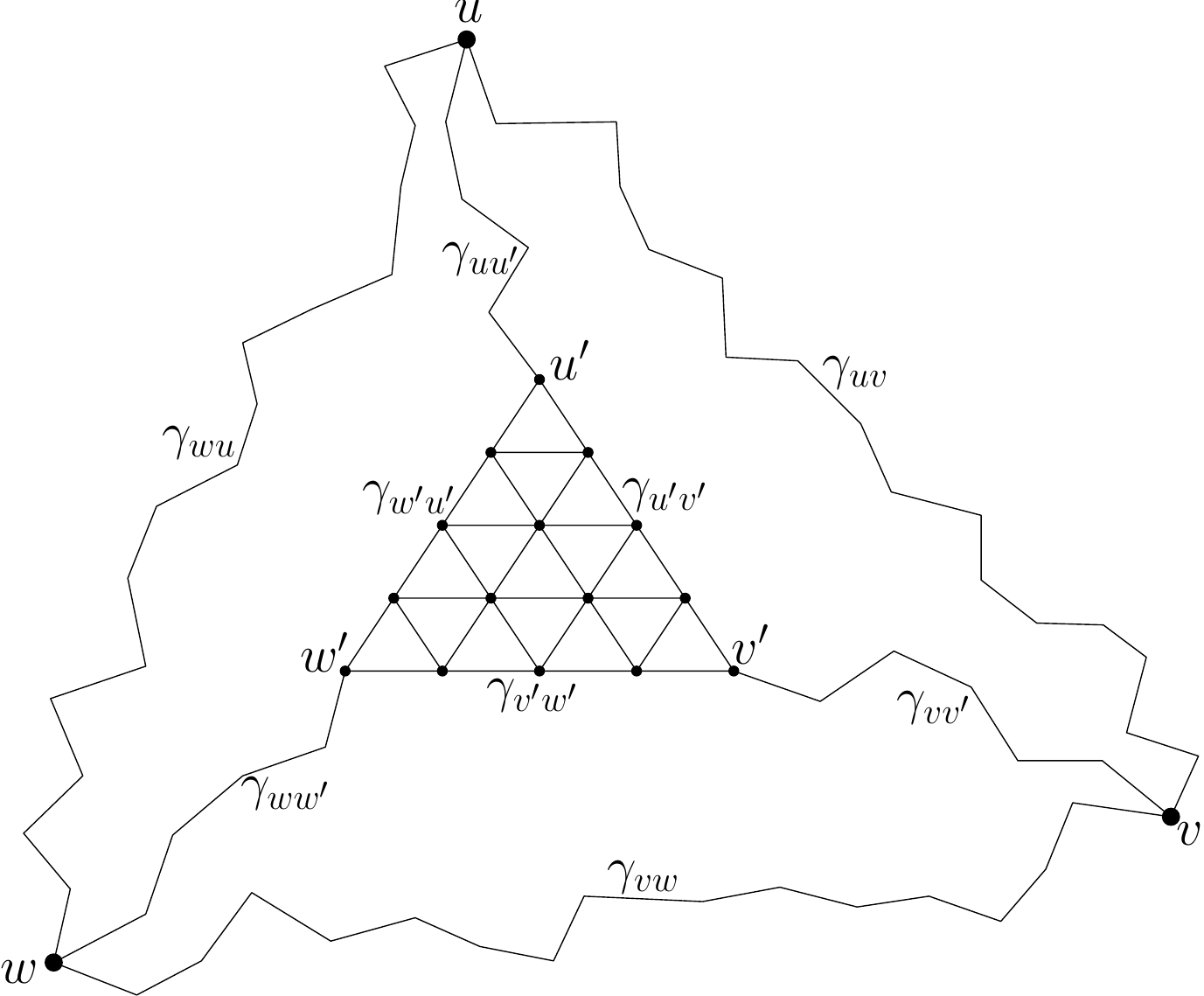}
\end{center}
\caption{Metric triangle.}\label{f:metrtr}
\end{figure}
By Lemma~\ref{l:metrsys}, there is a disc diagram $D_0 \to X$ for the loop $\gamma_{u'v'}\gamma_{v'w'}\gamma_{w'u'}$ with internal vertices degrees $6$ and boundary vertices of degree at most $4$.
By Lemma~\ref{l:bigon}, there is a disc diagram $D_1 \to X$ for the loop $\gamma_{uv}\circ \gamma_{vv'}\circ \gamma_{u'v'}^{-1}\circ \gamma_{uu'}^{-1}$ with internal vertices degrees $6$ and boundary vertices of degree at most $5$.
Analogously we obtain disc diagrams $D_2 \to X$ and $D_3 \to X$ for the bigons $vw$ and $wu$, respectively.
The required disc diagram $D \to X$ is then obtained as a combination of disc diagrams $D_i\to X$ for $D$ being 
the boundary union of discs $D_1,D_2,D_3,D_4$. It follows that the internal vertices in $D$ have degree 
at most $5+5+4=14$ and degrees of boundary vertices are bounded by $5$. Since the disc diagrams $D_i\to X$ are embeddings, we have that 
the map $D\to X$ is at most $4$--to--$1$.
\end{proof}
%\begin{thm}[{\cite[Proposition 6.6]{ChO}}]
%\label{l:fixpts}
%Let $G$ be a group acting by automorphisms on a weakly systolic complex $X$. Then the set Fix$_GX$ of points fixed by $G$ is contractible or empty. 
%\end{thm}

\begin{proof}[Proof of Proposition \ref{p:sap-wsys}]
	Having tight hexagons follows immediately from Proposition~\ref{p:tightsys} by decomposing a hexagon into four geodesic triangles. By Lemma~\ref{tightdiscsimplyrefinement} the Small Angle Property follows.
\end{proof}

\begin{proof}[Proof of Proposition \ref{p:fixpt-wsys}]
By \cite[Proposition 6.6]{ChO}, for any group $G$ acting by automorphisms on a weakly systolic complex
the set of points fixed by $G$ is contractible or empty.
\end{proof}

\begin{proof}[Proof of Proposition \ref{p:finint-wsys}]
Let $v,w$ be two vertices at distance $n$. We have to show that the interval between $v$ and $w$ is finite. By the vertex condition (V) from Definition~\ref{d:wsys}, and by the fact that all cliques in the $1$--skeleton are finite, there are finitely many vertices at distance $n-1$ from $v$ in the interval. Similarly we use the vertex condition inductively
to show that there are finitely many vertices in the interval at distance $n-k$ from $v$, for all $k\in {0,1,\ldots, n}$.
\end{proof}

%% file: sc.tex
\subsection{Combinatorial geometry of $C'(1/6)$ polygonal complexes}
\label{s:smallc}

The aim of this subsection is to prove Propositions~\ref{smallcancellationtight}, \ref{prop:SC_intervals}, and \ref{smallcancellationfixedpoint}, showing that polygonal complexes satisfying a metric small cancellation condition satisfy the combinatorial and algebraic conditions introduced in the previous sections.
As an immediate consequence of these results and Theorem A in Introduction, we obtain Theorem C.

\begin{prop}\label{smallcancellationtight}
A $C'(1/6)$ polygonal complex has tight hexagons. In particular, it satisfies the Small Angle Property.
\end{prop}
\begin{prop}\label{prop:SC_intervals}
	A $C'(1/6)$ polygonal complex has finite intervals.
\end{prop}
\begin{prop}\label{smallcancellationfixedpoint}
 For any action of a group on a  $C'(1/6)$ polygonal complex loops in fixed-point sets are contractible.
\end{prop}

We start by recalling some standard vocabulary about small cancellation theory.
% We show here that a $C'(1/6)$-polygonal complex in the sense of Wise \cite[Definition 2.4]{WiseSmallCancellation} satisfies the tight disc property, and hence the refinement property by Lemma \ref{tightdiscsimplyrefinement}. We refer to \cite{WiseSmallCancellation} for the usual background on small cancellation complexes, van Kampen diagrams, ladders, shells, etc.
% 
\begin{definition}[small cancellation]
%A \emph{combinatorial map} between CW-complexes is a map that restricts to an embedding on every open cell. A \emph{polygonal complex} is a $2$-dimensional CW-complex such that every $2$-cell has the structure of a polygon (i.e. $n$ vertices, $n$ edges, one open $2$-cell for some integer $n$) and such that every attaching map is combinatorial.
%
%Let $X$ be a polygonal complex. By a \textit{path of $\gamma$}, we mean an injective path in the $1$-skeleton of $Y$. The number of edges of a path $\gamma$ is denoted $|\gamma|$, and we refer to it as the \textit{length} of $\gamma$. 

Let $X$ be a polygonal complex. A \textit{piece} of $X$ is a path $\gamma$ such that there exist two polygons $R_1$ and $R_2$ of $X$ such that the map $\gamma \ra X$ factors as $\gamma \ra R_1 \ra X$ and $\gamma \ra R_2 \ra X$ but there does not exist a homeomorphism $\partial R_1 \ra \partial R_2$ that makes the following diagram commute: 
$$ \xymatrix{
     \gamma  \ar[d]_{} \ar[r]_{}^{} & \partial R_2 \ar[d]^{} \\
    \partial R_1 \ar[r]_{} \ar[ur]& X. \\
  }$$\\
By convention, edges of $X$ are also considered as pieces.

We say that $X$ is a $C'(1/6)$ \textit{polygonal complex} if every piece $\gamma$ of $X$ and every polygon $R$ of $X$ containing $\gamma$ satisfy the following relation:
$$|\gamma| < \frac{1}{6} \cdot |\partial R|.$$
\label{smallcancellationcomplex}
\end{definition}

In what follows, $X$ will denote   a given simply connected $C'(1/6)$-polygonal complex. 
%\begin{definition} Recall that a disc diagram $D$ over $X$ is a contractible planar polygonal complex, together with a combinatorial map  $D \ra X$ which embeds each polygon. The disc diagram $D$ over $X$ is called \textit{reduced} if no two distinct polygons of $D$ that share an edge are mapped to the same polygon of $X$.
%We denote by $\partial D$ the boundary of $D$.
%For a polygon $R$ of $D$, the intersection $\partial R \cap \partial D$ is called the \textit{outer component} of $R$ and denoted $\partial_o R$ (or the \textit{outer path}, when such an intersection is connected). The closure of $\partial R \setminus \partial_o R$ is called the \textit{inner component} of $R$ (or the \textit{inner path} when $\partial R \setminus \partial_o R$ is connected).
%
%A disc diagram is \textit{non-degenerate} if its boundary is homeomorphic to a circle, and \textit{degenerate} otherwise. An \textit{arc} of $D$ is a path of $D$ the interior vertices of which have valence $2$ and the extremities of which have valence at least $3$. Such an arc is called \textit{internal} if its interior is contained in $D \setminus \partial D$, and \textit{external} if it is contained in $\partial D$.
%\end{definition}
%
%Recall that, by a fundamental theorem of Lyndon--van Kampen, every loop of $X$ is the boundary of a reduced disc diagram (see for instance \cite{LyndonvanKampen}). 
We recall a fundamental combinatorial tool to study $C'(1/6)$ polygonal complexes.

\begin{definition}
Let $D$ be a planar contractible polygonal complex. 

A \textit{spur} of $D$ is an edge of $D$ with a vertex of valence $1$.

A \textit{shell} of $D$ is a polygon of $D$ such that $\partial R \cap \partial D$ is connected and whose inner path is a concatenation of at most $3$ internal arcs of $D$. 

The disc diagram $D$ is called a \textit{ladder} if it can be written as a union $D = c_1 \cup \ldots \cup c_n$, where the $c_i$ are distinct $1$- or $2$-cells such that:
\begin{itemize}
\item both $D \setminus c_1$ and $D \setminus c_n$ are connected,
\item the subspace $D\setminus c_i$ has exactly two connected components for $1 < i < n$,
\item if some $c_i$ is an edge, then no other $c_j$ contains it.
\end{itemize}
A path of the form $c_i \cap c_{i+1}$ for two consecutive $2$-cells of $D$ is called a \emph{rung}. For $1< k <n$ such that $c_k$ is a $2$-cell, the closure of a connected component of $c_k \setminus (c_{k-1} \cup c_{k+1})$ is called a \emph{rail} of $c_k$.
\end{definition}

\begin{thm}[Classification Theorem for disc diagrams in a $C'(1/6)$ polygonal complex \cite{McCammondWiseFansLadders}]
Let $D \ra X$ be a reduced disc diagram over $X$. Then one of the following holds:
\begin{itemize}
\item $D$ consists of a $0$-, $1$-, or $2$-cell,
\item $D$ is a ladder,
\item $D$ contains at least three shells or spurs.\qed
\end{itemize}
\label{classification}
\end{thm}

%\begin{cor}[{\cite[Lemma 13.2]{MaccammondWiseFansLadders}}]
%A reduced disc diagram over $X$ restricts to an embedding on every closed $2$-cell.
%\end{cor}

%\begin{lem}
%There are only finitely many combinatorial geodesics between two given vertices of a $C'(1/6)$ complex. More precisely, there exists a ladder containing every combinatorial geodesic between two given vertices. \qed
%\end{lem}
%
%\begin{cor}
% The family of combinatorial geodesics in a $C'(1/6)$ polygonal complex is a good system of paths. \qed
%\end{cor}

In order to show that $X$ has tight hexagons, it is necessary to have a finer understanding of the geodesic triangles of $X$. We will need the following:

\begin{definition}[geodesic ladder]
Let $D \ra X$ be a reduced disc diagram such that $D$ is a  non-singular ladder. We say that the disc diagram is a \textit{geodesic ladder} if its boundary path can be written as a union $\partial D = \gamma_1 \cup \gamma_2 \cup \gamma_1' \cup \gamma_2'$, where $\gamma_1', \gamma_2'$ are contained in the two shells of $D$ (if $D$ contains at least two cells, we require $\gamma_1', \gamma_2'$ to be in different shells), and $\gamma_1, \gamma_2$ map to geodesics of $X$.
\end{definition}

\begin{lem}
A reduced geodesic ladder embeds in $X$.
\label{geodesicladder}
\end{lem}

\begin{proof}
By contradiction, consider a reduced geodesic ladder $D$ that does not embed in $X$, and let $v,v'$ be two distinct vertices of $D$ that are sent to the same vertex $w$ of $X$. Let $\gamma$ be a path of $D$ from $v$ to $v'$ of the form $\gamma = A \cup P \cup B$, where $A, B$ are contained in a $2$-cell and $P$ is a subpath of $\gamma_1$ or $\gamma_2$. Then $\gamma$ maps to a loop of $X$, which we can assume to be injective, and let us choose a reduced non-singular disc diagram $D'\ra X$ with that loop as a boundary. 
First notice that $D'$ cannot be a single $2$-cell. Indeed, if that was the case, then $A$ and $B$ would be pieces, and $P$ being a geodesic contained in the boundary of $D'$, we would get 
$$|\partial R| = |A| + |P| + |B| < \frac{1}{6}|\partial R| + \frac{1}{2}|\partial R| + \frac{1}{6}|\partial R| < |\partial R|,$$
a contradiction. 
Thus, by the classification of diagrams, $D'$ contains at least two shells, and one of them does not contain $w$ in the interior of its outer path. Thus the shell $R$ \textit{lifts} to $D$, meaning we can form the new reduced disc diagram $D \cup R \ra X$. Notice that $R \cap D$ cannot contain a whole rail of a polygon $R'$ of $D$, for otherwise this path would be a piece, and since rungs also are pieces, the opposite rail of $R'$ would be of length strictly bigger than $\frac{1}{2} |\partial R'|$ by the $C'(1/6)$ condition, contradicting the fact that $D$ is a geodesic ladder. Thus $\partial_o R$ is the concatenation of at most two pieces, hence $\partial R$ is the concatenation of at most $5$ pieces since it is a shell of $D'$, a contradiction. 
\end{proof}

% \begin{prop}
% Let $X$ be a $C'(1/6)$ complex. Then $X$ has tight discs.
% \end{prop}

\begin{proof}[Proof of Proposition \ref{smallcancellationtight}]
 %alex: quote Dominik.
Let $\Delta(a,b,c)$ be an embedded geodesic triangle of $X$, and let $D \ra X$ be a (non-singular) reduced disc diagram with $\Delta$ as its boundary. By the classification of non-singular reduced disc diagrams filling a geodesic triangle  in a $C'(1/6)$ polygonal complex due to Strebel \cite[p.\ 261]{Strebel}, $D$ can be written as the union  of three geodesic ladders (the \textit{tails}) and at most three other $2$-cells (the \textit{core}). 

Since each tail is embedded by Lemma \ref{geodesicladder}, it follows that the map $D \ra X$ is at most 6-to-1. Furthermore,  it follows from the classification of such geodesic triangles that  each vertex has degree at most $3$. Thus $X$ has tight hexagons and, by Lemma \ref{tightdiscsimplyrefinement}, it satisfies the Small Angle Property.
\end{proof}

% \begin{lem}
% Let $G$ be a group acting by combinatorial isomorphisms on a $C'(1/6)$ complex $X$, and $H$ a subgroup of $G$. Then every loop contained in the fixed-point set $X^H$ is nullhomotopic.
% \end{lem}

\begin{proof}[Proof of Proposition \ref{prop:SC_intervals}] The fact that $C'(1/6)$ small cancellation polygonal complexes have finite intervals follows from \cite[Proposition 3.6]{GruberSisto}. The proof therein is given for Cayley graphs of classical $C'(1/6)$ small cancellation groups, but the same proof goes through for $C'(1/6)$ polygonal complexes.
%	We claim that, given two vertices of $C'(1/6)$ complex $X$, there exists a ladder that contains the interval between those vertices. Note that this claim will immediately proves the Proposition. Let $v, w$ be two vertices of $X$ and choose $\gamma$ a geodesic between them. For any other geodesic $\gamma'$ between $v, w$, Strebel's classification of bigons \cite{Strebel} implies that $\gamma$ and $\gamma'$ bound a geodesic ladder, so by Lemma \ref{geodesicladdereodesicladder} we will identify this disc diagram with its image in $X$. Let $\gamma_1, \gamma_2$ be two geodesics between $v$ and $w$. LLet $e$ be an edge of $\gamma$ and suppose that $e$ is contained in a polygon $P_1$ of the associated disc diagram $D_1$, and similarly $e$ is contained in a polygon $P_2$ of the associated disc diagram $D_2$. If $P_1$ were distinct from $P_2$, then $P:= P_1\cup P_2$ is the reunion of two polygons of a $C'(1/6)$ complex, hence is convex. What is more, $P$ contains an edge $e_1$ of $\gamma_1$
\end{proof}

\begin{proof}[Proof of Proposition \ref{smallcancellationfixedpoint}]
Let $G$ be a group acting by combinatorial isomorphisms on a $C'(1/6)$ complex $X$, and $H$ a subgroup of $G$. Let $\gamma$ be a loop in the $1$-skeleton of $X$ (which we can assume to be injective), which is pointwise fixed by $H$, and let $D \ra X$ be a reduced disc diagram with $\gamma$ as boundary. By the Classification Theorem \ref{classification},  there exists a polygon $R$ of $D$ whose outer path is connected and which is of length strictly more than $\frac{1}{2}|\partial R|$. In particular, such an outer path cannot be a piece by the small cancellation condition. As it is pointwise fixed by $H$, the whole of $\partial R$ is pointwise fixed by $H$, hence so is $R$. Thus $D \setminus R \ra X$ yields a disc diagram with strictly smaller area and whose boundary is pointwise fixed by $H$. The result now follows by induction on the number of polygons of $D$. 
\end{proof}

%% file: example.tex
%\section{New examples of hyperbolic groups}

\subsection{An example: Small cancellation over a graph of hyperbolic groups}
\label{s:scog}

As an application of our approach, we now prove Theorem D from Introduction. We consider a small cancellation quotient $\rquotient{G}{\ll\cR\gg }$ satisfying the assumptions of Theorem D.
%
%\begin{thm}\label{thm:hyperbolicity_small_cancellation}
% Let $G(\Gamma)$ be a finite graph of groups over a simplicial graph $\Gamma$ satisfying the following: 
% \begin{itemize}
%  \item every vertex group is hyperbolic,
%  \item every edge group embeds as a quasiconvex malnormal subgroup in the associated vertex groups,
%  \item for every vertex $v$ of $\Gamma$, the family of adjacent edge groups $(G_e)_{v \in e}$ is almost malnormal.
%   \end{itemize}
%Let $G$ be the fundamental group of $G(\Gamma)$. Let $\cR$ be a finite set of relators satisfying the classical $C'(1/6)$--small cancellation over $G(\Gamma)$. Then the quotient group $\rquotient{G}{\ll\cR\gg }$ is hyperbolic and the quotient map $G \ra \rquotient{G}{\ll\cR\gg }$  embeds every local group of $G(\Gamma)$ as a quasiconvex subgroup.  
%\end{thm}
In order to do so, we use an action of   $\rquotient{G}{\ll\cR\gg }$  to which Theorem \ref{thm:combination} applies. In the same way that classical $C'(1/6)$ small cancellation groups act geometrically on $C'(1/6)$ polygonal complexes, $C'(1/6)$ small cancellation groups obtained by killing off finitely many relations act cocompactly on $C'(1/6)$ polygonal complexes. Such constructions  of actions are well known to experts and can be found detailed in several places, for instance in  \cite{MartinSteenbockFreeProduct} in the particular case of small cancellation over free products, or in \cite{MartinSmallCancellationClassifying} for $C''(1/6)$ small cancellation over graphs of groups. In \cite{MartinSmallCancellationClassifying}, the strategy  can be thought of as starting from the action on a small cancellation complex associated to the group, and identifying certain subspaces to obtain an action on a CAT(0) space, following an idea of Gromov \cite{GromovCATkappa}. In order to obtain such a CAT(0) space, the constructions from \cite{MartinSmallCancellationClassifying} were made in the  $C''(1/6)$ setting, but the construction of the action on a  polygonal complex would work in the weaker $C'(1/6)$ setting we are considering here: this is the construction of the space denoted $X$ in  \cite[Definition 6.3]{MartinSmallCancellationClassifying}.  In particular, the following holds:

\begin{thm}\label{thm:MartinSteenbock}
The small cancellation quotient $\rquotient{G}{\ll\cR\gg }$ acts cocompactly on a hyperbolic $C'(1/6)$ small cancellation complex $X$ such that: 
\begin{itemize}
	\item the stabiliser of a vertex of $X$ is isomorphic to a vertex group of $G(\Gamma)$,
	\item the action is without inversion on the $1$-skeleton of $X$, and the stabiliser of an edge  of $X$ is isomorphic to an edge group of $G(\Gamma)$,
	\item global stabilisers of polygons of $X$ are finite,
	\item for every vertex $v$ of $X$, there exists a vertex $v'$ of $\Gamma$ such that for every edge $e$ of $X$ containing $v$, there exists an edge $e'$ of $\Gamma$ containing $v'$ such that the inclusion $G_e \hra G_v$ is conjugated to the local morphism $G_{e'} \hra G_{v'}$.\qed
\end{itemize}	
\end{thm}

\begin{proof}[Proof of Theorem D]
%Since the aforementioned complex $X$ is a $C'(1/6)$ polygonal complex, it has the small angles property by Proposition~\ref{smallcancellationtight}, finite intervals by Proposition~\ref{prop:SC_intervals}, and satisfies the fixed point property by Proposition~\ref. Moreover, s
Since for every vertex $v$ of $\Gamma$, the family of subgroups $G_e$, where $e$ is an edge of $\Gamma$ containing $v$, is almost malnormal, it follows from Theorem \ref{thm:MartinSteenbock} that the action of $\rquotient{G}{\ll\cR\gg }$ on $X$ is weakly acylindrical. 
Since the aforementioned complex $X$ is a $C'(1/6)$ polygonal complex, the result now follows from Theorem \ref{thm:MartinSteenbock}, Theorem C and Theorem A.
\end{proof}

%% file: appendices.tex
\begin{appendices}

\section{The topology of $\partial G$ and the dynamics of the associated action}

In this appendix, we prove Theorem \ref{thm:combination} by adapting the proofs of \cite{MartinBoundaries} to our combinatorial setting. The proofs presented in \cite{MartinBoundaries} extend to this case with little changes, and we prove here the results that need slight modifications,   following very closely the structure of the original proofs. Every time a result parallels a result of \cite{MartinBoundaries} we give a reference to that analogous result and, when appropriate, we explain how the new combinatorial tools developed in this article allow us to adapt the proofs to this new setting.

In this appendix, we consider a group $G$ acting on a complex $X$ satisfying the hypotheses of Theorem A, and we denote by $\delta$ the hyperbolicity constant of $X$. We start by recalling the changes between \cite{MartinBoundaries} and the present articles: 
\begin{itemize}
	\item In 	\cite[Section 2.1]{MartinBoundaries}, the compactification of $G$ being  considered was $EG \sqcup \partial G,$ where $EG$ was a classifying space for proper actions of $G$. Instead, we consider here the compactification  $\cV G \sqcup \partial G.$ Having points of $\cV G$ being isolated points actually makes some of the proofs  (separation of neighbourhoods, convergence results, etc.) much easier and shorter.
  \item In \cite[Section 6.1]{MartinBoundaries}, the topology on the compactification was defined using the behaviour of CAT(0) geodesics. Here instead, we consider combinatorial geodesics. A difference is that there can be several geodesics between two different vertices. However, if $\gamma, \gamma'$ are two geodesics between vertices $v$ and $w$, the Small Angle Property ensures that there is a path of simplices around $v$ of length bounded by a universal constant between the first edges of $\gamma$ and $\gamma'$. Thus, considering intervals instead of CAT(0) geodesics poses no significant problem in adapating the proof of \cite{MartinBoundaries}.
	\item  In \cite[Definition 5.1]{MartinBoundaries}, the first author introduced the notion of exit simplex $\sigma_{\xi, \varepsilon}(\gamma)$, that is, the first simplex met by a CAT(0) geodesic $\gamma$ when leaving the $\varepsilon$-neighbourhood of the domain $D(\xi)$. Instead, adopting a purely combinatorial approach here, we just have to consider combinatorial geodesics, and in particular we only consider exit edges $e_\xi(\gamma)$. Again, this simplification greatly shortens some of the proofs.
	\item In \cite[Definition 6.5]{MartinBoundaries}, neighbourhoods of points of $\partial_{Stab}G$, denoted $V_{\cU, \varepsilon}(\xi)$, were defined by considering the way CAT(0) geodesics exit the $\varepsilon$-neighbourhood of the domain $D(\xi)$. Here again, we only consider combinatorial geodesics, and consequently neighbourhoods are defined in terms of the way combinatorial geodesics exit the domain. Such neighbourhoods are denoted  $V_{\cU}(\xi)$, and again, this approach simplifies some of the proofs. Analogously, the cones $Cone_{\cU, \varepsilon}(\xi)$ considered in \cite{MartinBoundaries} are replaced here by their combinatorial counterpart $Cone_{\cU}(\xi)$.
\end{itemize}

In what follows, most of the proofs from \cite{MartinBoundaries} translate almost immediately by replacing the CAT(0) notions by their combinatorial counterpart, as explained above. 

 We start by proving that $\partial G$ is a compact metrisable space. It follows closely the structure of \cite[Section 7]{MartinBoundaries}, and we refer to the corresponding result whenever possible.

\subsection{Nestings}

\begin{definition}[nested]
 Let $\xi \in \partial_{Stab} G$, $v$ be a vertex of $D(\xi)$, and $U$ be a neighbourhood of $\xi$ in $\cC G_v$. We say that a subneighbourhood $V \subset U$ containing $\xi$ is \textit{nested} in $U$ if its closure is contained in $U$ and for every simplex $\sigma$ of $\mbox{st}(v)$ (star of $v$) not contained in $D(\xi)$, we have
$$ \cC G_\sigma\cap V \neq \varnothing \Rightarrow \cC G_\sigma \subset U.$$
\end{definition}

\begin{rmk}
For every point $\xi$ of a fibre $\partial G_v$ and every neighbourhood $U$ of $\xi$ in $\partial G_v$, there exists a subneighbourhood $V$ of $U$ containing $\xi$ which is nested in $U$ \cite[Lemma 4.10]{MartinBoundaries}
\end{rmk}

\begin{definition}[Definition 4.13 of \cite{MartinBoundaries}]
 Let $\xi \in \partial_{Stab} G$, together with two $\xi$-families $\cU,  \cU'$. We say that $\cU'$ is \textit{nested} in $\cU$ if for every vertex $v$ of $D(\xi)$, $U_v'$ is nested in $U_v$. Furthermore we say that $\cU'$ is $n$-nested in $\cU$ if there exist $\xi$-families
$$\cU' = \cU^{[0]} \subset \ldots \subset \cU^{[n]} = \cU$$
with $\cU^{[i]}$ nested in $\cU^{[i+1]}$ for every $i=0,\ldots,n-1$.
\label{defnested}
\end{definition}

% \begin{definition}[$\xi$-family]
%  Let $\xi \in \partial_{Stab} G$. A collection $\cU$ of open sets $\left\lbrace U_v, v \in V(\xi) \right\rbrace $ is called a $\xi$-\textit{family} if for every pair of vertices $v,v'$ of $X$ that are joined by an edge $e$ and every $x \in \overline{EG_e}$,
% $$\phi_{v,e}(x) \in U_v \Leftrightarrow \phi_{v',e}(x) \in U_{v'}.$$
% \label{defxifamily}
% \end{definition}
% 
% \begin{prop}[Proposition ??? of \cite{MartinBoundaries}]
%  Let $\xi \in \partial_{Stab} G$. For every vertex $v$ of $D(\xi)$, let $U_v$ be a neighbourhood of $\xi$ in $\overline{EG_v}$. Then there exists a $\xi$-family $\cU'$ such that $U_v' \subset U_v$ for every vertex $v$ of $D(\xi)$.\qed
% \label{xifamilies}
% \end{prop}

Recall that domains contain at most $d_{\mathrm{max}}$ by Definition \ref{def:dmax}. Recall also that, following Definition \ref{def:small_angle_property}, for every integer $n \geq 0$ there exists a constant $r(n)$  so that the Small Angle Property holds for subcomplexes of $X$ containing at most $n$ edges.

\begin{definition}[refined family]
Let $\xi$ be a point of $\partial_{Stab} G$, $\cU$ a $\xi$-family. A $\xi$-family $\cU'$ which is $r(d_{max})$-nested in $\cU$ is said to be \textit{refined} in $\cU$. Furthermore we say that $\cU'$ is $n$-refined in $\cU$ if there exist $\xi$-families
$$\cU' = \cU^{[0]} \subset \ldots \subset \cU^{[n]} = \cU$$
with $\cU^{[i]}$ refined in $\cU^{[i+1]}$ for every $i=0,\ldots,n-1$.
\end{definition}

%An immediate corollary is the following: 
%
%\begin{lem}
%	Let $\xi$ be a point of $\partial_{Stab} G$, $\cU$ a $\xi$-family. 
%	\end{lem}

\begin{rmk}
 In \cite{MartinBoundaries}, a refined $\xi$-family was called $d_{max}$-refined. The notation used here is slightly less cumbersome. 
\end{rmk}

\begin{definition}[$\xi$-family not seeing some subspace]
 Let $\xi$ be a point of $\bds$, $\cU$  a $\xi$-family and $K$ a set of generalised vertices of $\overline{X}$. We say that $\cU$ \emph{does not see} $K$ if the following holds: for every geodesic $\gamma$ between a vertex $v \in D(\xi)$ and a point $z \in K$, the unique edge $e$ of $\gamma$ contained in $  N(D(\xi))\setminus D(\xi)$ is such that, if we denote $w:= e \cap D(\xi)$, we have  
 $$ \mathcal{C}G_e \cap U_w = \varnothing.$$
% $$(K \setminus D(\xi)) \cap \mbox{Cone}_{\cU}(\xi) = \varnothing.$$
\end{definition}

\begin{lem}
Let $\xi$ be a point of $\bds$ and $K$ a set of generalised vertices of $\overline{X}$. There exists a $\xi$-family which does not see $K$ in the following two cases:
\begin{itemize}
 \item $K$ is a finite set of vertices of $X$,
 \item $K$ consists of a single point of $\partial X$.
\end{itemize}
\label{doesnotsee}
\end{lem}

\begin{proof}
First consider the case where $K$ is a finite set of vertices of $X$. Since combinatorial intervals are finite by assumption and since domains are finite subcomplexes by Proposition \ref{finitedomain}, there are only finitely many geodesics between $D(\xi)$ and $K$. For every vertex $v$ of $D(\xi)$, let $\cE_v$ be the set of exit edges of geodesics from $D(\xi)$ to a point of $K$ which leave the domain $D(\xi)$ at the vertex $v$. We can thus find a neighbourhood $U_v$ of $\xi$ in $\CG_v$ which is disjoint from every $\CG_e \subset \CG_v$, $e \in \cE_v$. Now any $\xi$-family contained in the family of $U_v, v \in D(\xi)$ works, by construction of a cone.

Consider now the case of an element $\eta \in \partial X$. Let $N \geqslant 0$ be an integer such that $D(\xi)$ is contained in the $N$-ball around $v_0$. Choose a based geodesic $\gamma_{\eta}$ from $v_0$ to $\eta$ and set $x:= \gamma_\eta(N+\delta +1)$. Choose a $\xi$-family $\cU$ which does not see $x$ and a $\xi$-family $\cU'$ that is %$\delta$-
refined in $\cU$.  In particular, $\cU'$ does not see the $\delta$-ball around $x$. By definition of $\delta$, every geodesic from $v_0$ to $\eta$ meets the ball $B(x,\delta)$, it thus follows that $\cU'$ does not see $\eta$.
\end{proof}

\begin{conv}
	From now on, we will assume that $\xi$-families do not see the basepoint $v_0$.
	\end{conv}

\subsection{The geometric toolbox}
\label{sec:toolbox}

An important part of \cite{MartinBoundaries} was to develop sufficiently fine tools to understand the topology of the space $\cC G$. The main results used are contained in the following `toolbox'. While the CAT(0) geometry of the space was used in a few instances (in which case we will explain how to adapt the proof), these lemmata really represent the main geometric tools, and the proofs carry over in our new combinatorial effortlessly once the appropriate analogues are available. 

In this paper, an important geometric tool is the Small Angle Property. This is a combinatorial analogue of a result from CAT(0) geometry proved in \cite{MartinBoundaries}, namely the Short Path of Simplices Lemma \cite[Lemma 3.7]{MartinBoundaries}. This lemma was crucial in proving the main geometric tools from \cite{MartinBoundaries}, and likewise the Small Angle Property plays a key role in proving the combinatorial analogues of these lemmata, which we now introduce. 

\begin{lem}[Geodesic Reattachment Lemma, compare with {\cite[Lemma 5.8]{MartinBoundaries}}]
	Let $\xi \in \partial_{Stab}G$, $\cV$ a $\xi$-family that does not see $v_0$, $\cU$ a $\xi$-family that is   refined in $\cV$, and $x \in X \setminus D(\xi)$. Suppose that there exists an edge $e$ of $N(D(\xi)) \setminus D(\xi) $ contained in a geodesic $\gamma$ from $x$ to a vertex of $ D(\xi)$,  such that for the  vertex $v:= e \cap D(\xi)$, we have $\cV G_e \subset U_v$. Then any geodesic from $v_0$ to $x$ meets  $D(\xi)$, and we have  $x\in  Cone_\cV(\xi)$.
	\label{geodesicreattachment}
\end{lem}	

\begin{proof}
	By contradiction, let $\gamma'$ be a geodesic from $v_0$ to $x$ that does not meet $D(\xi)$, and let $\gamma_0$ be a geodesic from $v_0$ to $v$. Let us denote by $e_0$ the unique edge of $\gamma_0$ contained in $N(D(\xi)) \setminus D(\xi)$, and let us denote $w_0:= D(\xi)\cap e_0$. Since $D(\xi)$ contains at most $d_\mathrm{max}$ simplices, it follows from the Small Angle Property that there exists a path of simplices of length at most $r(d_\mathrm{max})$ between $e$ and $e_0$. Since $\cU$ is $r(d_\mathrm{max})$-nested in $\cV$, it follows that $\partial G_{e_0} \subset V_{w_0}$, which contradicts the fact that $\cV$ does not see $v_0$.  
	
	Thus, every geodesic from $v_0$ to $x$ meets $D(\xi)$. Moreover, by the Small Angle Property  there exists a path of simplices of length at most $r(d_\mathrm{max})$ between $e$ and $e_\xi(\gamma)$. Since $\cV G_e \subset U_v$ and $\cU$ is refined in $\cV$, it follows that $x\in  Cone_\cV(\xi)$.
	\end{proof}

\begin{lem}[Refinement Lemma, compare with {\cite[Lemma 5.10]{MartinBoundaries}}]
	Let $\xi \in \partial_{Stab}G$, $\cV$ a $\xi$-family that does not see $v_0$, $\cU$ a $\xi$-family that is   $n$-refined in $\cV$. Then the following holds: 
	
	Let $\tau$ be a  combinatorial path being a concatenation of at most $n$ geodesics in $X \setminus D(\xi)$. Assume that there exists a vertex $v$ of $\tau$ such that there exists a geodesic $\gamma$ from $v_0$ to $v$ meeting $D(\xi)$ and such that $e_\xi(\gamma)\subset U_{v_\xi(\gamma)}$. Then $\tau \subset Cone_\cV(\xi)$.
	\label{refinement}
	\end{lem}

\begin{proof}
	This is a straightforward induction on $n$, using the Geodesic Reattachment Lemma \ref{geodesicreattachment} to each of the $n$ geodesic segments composing $\tau$, first to show that any geodesic segment from $v_0$ to a point of $\tau$ meets $D(\xi)$, and then to control the way these geodesics exit $D(\xi)$.
\end{proof}

\subsection{Basis of neighbourhoods}

We now prove the following: 

\begin{thm}[Basis of neighbourhoods, see {\cite[Theorem 6.17]{MartinBoundaries}}]
The family $\cB_{\cC G}$ is a basis of neighbourhoods for the topology of $\cC G$.
\end{thm} 

 In order to do so, we first prove the following:

\begin{prop}[Filtration, compare with {\cite[Filtration Lemma in Section 6.2]{MartinBoundaries}}] Let $z$ be a point of $\cC G$ and $U\in \cB_{\cC G}(z)$ a neighbourhood of $z$. Then there exists a subneighbourhood $U' \in \cB_{\cC G}(z)$ of $U$ such that every element 
$z' \in U'$ admits a neighbourhood $U'' \in \cB_{\cC G}(z')$ which is contained in $U$.
\end{prop}

The proof is very similar to that of \cite{MartinBoundaries}. As is the case there, it splits in many cases.

\begin{lem}
Let $z, z'$ be  points of $\cV G$,  and $U$ be an element of $\cB_{\cC G}(z)$ such that $z' \in U$. Then there exists a  neighbourhood $U' \in \cB_{\cC G}(z')$ such that $U' \subset U$.\qed
\label{filtrationx}
\end{lem}

 \begin{proof}
 	Take $U'= U$.
 \end{proof}
 
\begin{lem}
Let $\eta, \eta'$ be two points of $\partial X$ and $U$ be an element of $\cB_{\partial X}(\eta)$ such that $\eta' \in V_U(\eta)$. Then there exists a  neighbourhood $U' \subset \cB_{\partial X}(\eta')$ such that $V_{U'}(\eta') \subset V_U(\eta)$.
\label{filtrationetaeta}
\end{lem}

\begin{proof}
Take any neighbourhood $U' \in \cB_{\partial X}(\eta')$ contained in $U$.
\end{proof}

\begin{lem}
Let $\eta$ be a point of $\partial X$, $\xi$ be a point of $\partial_{Stab} G$ and $U$ be an element of $\cB_{\partial X}(\eta)$ such that $\xi \in V_U(\eta)$. Then there exists a  $\xi$-family $\cU$ such that $V_\cU(\xi) \subset V_U(\eta)$. 
%$U' \subset \cB_{\partial X}(\eta')$ such that $V_{U'}(\eta') \subset V_U(\eta)$.
\label{filtrationetaxi}
\end{lem}

\begin{proof}
Let $\cU$ be any $\xi$-family and $z$ be a point of $V_\cU(\xi)$. Let $x$ be a point of $p(z)$. Any geodesic from $v_0$ to $x$ goes through $D(\xi)$ by definition, and $D(\xi) \subset U$. By definition of $\cB_{\partial X}(\eta)$, it then follows that $x \in U$, thus $z \in V_U(\eta)$.
\end{proof}

\begin{lem}
Let $\xi$ be a point of $\partial_{Stab} G$ and $\cU$ a $\xi$-family. Then there exists a $\xi$-family $\cU'$ such that for every point $\eta \in \partial X$ with $\eta \in V_{\cU'}(\xi)$, there exists a neighbourhood $U'' \in \cB_{\partial X}(\eta)$ such that $V_{U''}(\eta) \subset V_\cU(\xi)$. 
%$U' \subset \cB_{\partial X}(\eta')$ such that $V_{U'}(\eta') \subset V_U(\eta)$.
\label{filtrationxieta}
\end{lem}

\begin{proof}
Let $\cU'$ be a $\xi$-family which is %$\delta$-
refined in $\cU$ and let $\eta \in \partial X$ with $\eta \in V_{\cU'}(\xi)$. Since $D(\xi)$ is finite by Proposition \ref{finitedomain}, let $N \geqslant 0$ be an integer such that $D(\xi)$ is contained in the $N$-ball around $v_0$. Let $\gamma_\eta$ be a geodesic from $v_0$ to $\eta$ and $x:= \gamma_\eta(N+\delta +1)$. Let 
$$U'' = \{ y \in \overline{X} \mbox{ such that } \langle y, \eta \rangle_{v_0} > N+ 2\delta +1 \} \subset \cB_{\overline{X}}(\eta)$$
and let $z$ be an element of $V_{U''}(\eta)$. Let $y$ be a point of $p(z)$, $\gamma_y$ be a geodesic from $v_0$ to $y$ and $y':= \gamma_y(N+\delta +1)$. By definition of $U''$ and since $X$ is $\delta$-hyperbolic, we have $d(x,y')\leqslant 2\delta$. Now any geodesic path from $x$ to $y'$ misses the $N$-ball around $v_0$ by construction, hence misses $D(\xi)$. The  Refinement Lemma \ref{refinement} thus implies that $y'$, hence $y$, is in $\mbox{Cone}_\cU(\xi)$, and thus $z \in V_\cU(\xi)$.
\end{proof}

\begin{lem}
Let $\xi$ be a point of $\partial_{Stab} G$, $\cU$ a $\xi$-family and $\xi'$ a point of $V_\cU(\xi')$. Then there exists a $\xi$-family $\cU'$ such that $V_{\cU'}(\eta') \subset V_\cU(\xi)$.
\label{filtrationxixi}
\end{lem}

\begin{proof}
Let $\cU'$ be a $\xi'$-family that does not see $D(\xi)$ and let us prove that $V_{\cU'}(\xi') \subset V_\cU(\xi)$. Let $z$ be a point of $V_{\cU'}(\xi')$, $y$ be a point of $p(z)$ and $\gamma_y$ be a geodesic from $v_0$ to $y$. By convexity of domains (Proposition \ref{finitedomain}) and construction of $\cU'$, $\gamma_y$ cannot meet $D(\xi)$ after leaving $D(\xi')$. Thus, one of the following situations occurs:
\begin{itemize}
\item If $\gamma_y$ meets $D(\xi')$ after leaving $D(\xi)$, then since $D(\xi') \setminus D(\xi) \subset \mbox{Cone}_\cU(\xi)$, it follows that $y \in  \mbox{Cone}_\cU(\xi)$. 
\item If $\gamma_y$ leaves $D(\xi)$ and $D(\xi')$ at the same vertex $v$, then $\gamma_y$ exits $D(\xi)$ in the direction of $U_v' \subset U_v$, hence  $y \in \mbox{Cone}_\cU(\xi)$.
\item If $y \in D(\xi) \cap D(\xi')$, then $z \in U_v' \subset U_v$.
\end{itemize}
In every situation, it follows that $z \in V_\cU(\xi)$.
\end{proof}

\begin{thm}[see {\cite[Theorem 6.17]{MartinBoundaries}}]
$\cB_{\cC G}$ is a basis for the topology of $\cC G$. This turns $\cC G$ into a second countable space.
\label{basis}
\end{thm}

\begin{proof}
By Lemmas \ref{filtrationx}, \ref{filtrationetaeta}, \ref{filtrationetaxi}, \ref{filtrationxieta} and \ref{filtrationxixi}, $\cB_{\cC G}$ satisfies the Filtration Property, hence is a basis of neighbourhoods. Let us now prove that the topology it defines is second countable. 

Since each $G_v$ is hyperbolic (and in particular finitely generated) and the action of $G$ on $X$ is cocompact, it follows that the simplicial complex $X$ has countably many cells, so the family of neighbourhoods $(V_n(x))_{n \geqslant 0, x \in \cV(X)}$  is countable.

A neighbourhood of a point $\xi$ of $\partial_{Stab} G$ is completely determined by its domain and the associated $\xi$-family. Domain are finite subcomplexes of $X$ by Proposition \ref{finitedomain}, so there are at most countably many of them. Furthermore, for every vertex $v$ of $X$, the space $\overline{EG}_v$ has a countable basis of neighbourhoods. From this it is clear that we can define a countable family of open neighbourhoods containing a basis of neighbourhoods of every point of $\partial_{Stab} G$.

Finally, there are countably many finite subsets of $\cV G$.
\end{proof}

\subsection{Induced topologies}

We have the following result, which is essentially \cite[Proposition 6.19]{MartinBoundaries}. As the proof of the aforementioned proposition does not use any CAT(0) geometry, it carries over to this combinatorial framework without any essential change.

\begin{prop}[compare with  {\cite[Proposition 6.19]{MartinBoundaries}}]
 The topology of $\CG$ induces the natural topologies on $\VG$, $\partial X$ and $\cC G_v$ for every vertex $v$ of $X$.\qed
\label{inducedtopology}
\end{prop}

\subsection{The $T_0$-condition}

We now prove that $\cC G$ is a $T_0$-space, following the proof of the analogous Proposition \cite[7.1]{MartinBoundaries}. Recall that this means that for every pair of distinct points of $\cC G$, there exists an open set that contains exactly one of them. We split the proof of the $T_0$ condition into different cases.  

\begin{lem}
Let $z \in \cC G$ be and $x \in \cV G$ be two distinct points. Then $z,x$ admit disjoint neighbourhoods.
\end{lem}

\begin{proof}
	Immediate as points of $\cV G$ are isolated by construction.
	\end{proof}

\begin{lem}
Let $\eta, \eta'$ be two distinct points of $\partial X \subset \partial G$. Then $\eta, \eta'$ admit disjoint neighbourhoods.
\end{lem}

\begin{proof}
The space $\overline{X} = X \cup \partial X$ is metrisable, hence Hausdorff. Disjoint neighbourhoods of $\eta, \eta'$ in $\overline{X}$ yields disjoint neighbourhoods of $\eta, \eta'$ in $\partial G$.
\end{proof}

\begin{lem}
Let $\eta \in \partial X$ and $\xi \in \partial_{Stab} G$. Then there exists a neighbourhood of $\eta$ which does not contain $\xi$.
\end{lem}

\begin{proof}
By Lemma \ref{doesnotsee}, choose a $\xi$-family $\cU$ that misses $\eta$. This defines a neighbourhood $V_\cU(\xi)$ not containing $\eta$.
\end{proof}

\begin{lem}
Let $\xi, \xi'$ be two distinct points of $\partial_{Stab} G$. Then they admit disjoint neighbourhoods.
\end{lem}

\begin{proof}
It is enough to choose a $\xi$-family $\cU$ that does not see $\xi' $ and a $\xi'$-family $\cU'$ that does not see $\xi$, and such that $\cU$ and $\cU'$ are disjoint. Such families exist by Lemma \ref{doesnotsee}.
\end{proof}

\begin{cor}[see {\cite[Proposition 7.1]{MartinBoundaries}}]
The space $\cC G$ satisfies the $T_0$ condition. \qed
\label{T0}
\end{cor}

\subsection{Regularity}

We now prove that $\cC G$ is regular, following the proof of \cite[Proposition 7.8]{MartinBoundaries}:

\begin{prop}[see {\cite[Proposition 7.8]{MartinBoundaries}}]
The space $\cC G$ is regular, that is, for every point $z \in \cC G$ and every neighbourhood  $U \in \cB_{\cC G}(z)$, there exists a subneighbourhood $U' \in \cB_{\cC G}(z)$ such that every point of $\cC G \setminus U$ admits a neighbourhood that is disjoint from $U'$.
\label{regular}
\end{prop}

We split the proof in three cases. 

\begin{lem}
Let $x \in \cV G$ and $U$ a neighbourhood of $x $ in $\cC G$. Then there exists a subneighbourhood $U'$ of $x$ such that every point of $\cC G \setminus U$ admits a neighbourhood disjoint from $U'$. 
\label{regular0}
\end{lem}

\begin{proof} Take $U'$ to be the neighbourhood consisting of the single point $x$. Now every point of $\cC G$ distinct from $x$ admits a neighbourhood that does not contain $x$: This is obvious for points of $\cV G$. For a point $\eta \in \partial X$, a neighbourhood $U''$ of $\eta$ in $X \cup \partial X$ not containing $p(x)$ yields a neighbourhood $V_{U''}(\eta)$ not containing $x$. For a point $\xi \in \partial_{Stab}G$, a $\xi$-family $\cU$ that does not see $p(x)$ yields a neighbourhood $V_{\cU}(\xi)$.
\end{proof}

\begin{lem}
Let $\eta \in \partial X$ and $U$ a neighbourhood of $\eta $ in $\overline{X}$. Then there exists a subneighbourhood $U'$ such that every point of $\cC G \setminus V_U(\eta)$ admits a neighbourhood disjoint from $V_{U'}(\eta)$.
\label{regular1}
\end{lem}

\begin{proof}
By Proposition \ref{finitedomain} there exists a constant $A$ bigger than the diameters of all domains.
Since $X$ is hyperbolic, there exists a subneighbourhood $U'$ of $U$ such that the distance between $X \setminus U$ and $ U' \cap X$ is strictly greater $A + \delta $, where 
%$A$ is the acylindricity constant and 
$\delta$ is the hyperbolicity constant. 
Since $\overline{X} $ is metrisable, hence regular, we can further assume that the closure of $U'$ is contained in $U$. Finally, we can assume that there exists an integer $N \geqslant 0$ such that 
$$ U' = \{ y \in \overline{X} \mbox{ such that } \langle y, \eta \rangle_{v_0} > N \}.$$

 For a point $x \in \cV G \setminus V_U(\eta)$, we just take the neighbourhood consisting of $x$ itself, which yields a neighbourhood of $x$ disjoint from $V_{U'}(\eta)$.

Let $\eta' \in \partial X \setminus U$. We have $\eta' \notin \overline{U'}$, so by regularity of $\overline{X}$ there exists a neighbourhood $U''$ of $\eta'$ in $\partial X$ disjoint from $U'$. This yields a neighbourhood $V_{U''}(\eta')$ disjoint from $V_{U'}(\eta)$.

Let $\xi \in \partial_{Stab}G \setminus V_U(\eta)$. Let $\cU$ be a $\xi$-family that does not see $\gamma_\eta(N)$ and $\cU'$ a $\xi$-family that is %$2d_\mathrm{max}$-
$2$--refined in $\cU$. %(where the constant $d_\mathrm{max}$ was defined in Definition \ref{???}).
Since $\xi \notin V_U(\eta)$, we have $D(\xi) \nsubseteq U$, hence $D(\xi) \cap U' = \varnothing$ since domains have a diameter bounded above by $A$. To prove that $V_{U'}(\eta)$ and $V_{\cU'}(\xi)$ are disjoint, it is enough to prove that $\cU'$ does not see $U'$. Let $x \in U'$, choose a geodesic $\gamma_x$ from $v_0$ to $x$. By definition of $U'$, we have $d(\gamma_x(N), \gamma_\eta(N)) \leqslant \delta$ and any geodesic from $\gamma_x(N)$ to $\gamma_\eta(N)$ is disjoint from $D(\xi)$ by construction of $U'$. Furthermore, the portion of $\gamma_x$ between $\gamma_x(N) $ and $x$ is also disjoint from $D(\xi)$ by construction. Thus, since $\cU'$ is %$2d_\mathrm{max}$-
$2$--refined in $\cU$, it follows from  the Refinement Property \ref{refinement} applied to $[\gamma_\eta(N), \gamma_x(N)] \cup [\gamma_x(N),x]$ that $\cU'$ does not see $x$.
\end{proof}

\begin{lem}
Let $\xi \in \partial_{Stab} G$ and $\cU$ a $\xi$-family. Then there exists a sub-$\xi$-family $\cU'$ such that every point of $\partial G \setminus V_\cU(\xi)$ admits a neighbourhood disjoint from $V_{\cU'}(\xi)$.
\label{regular2}
\end{lem}

\begin{proof}
Let $\cU'$ be a $\xi$-family that is %max$(d_{max}, (\lambda+1)\delta + \varepsilon)$-
refined in $\cU$.

For a point $x \in \cV G \setminus V_U(\eta)$, we just take the neighbourhood consisting of $x$ itself, which yields a neighbourhood of $x$ disjoint from $V_{\cU'}(\xi)$. 

Let $\eta \in \partial X \setminus V_\cU(\xi)$ and choose a geodesic $\gamma_\eta$ from $v_0$ to $\eta$. Let $N \geqslant 0$ be an integer such that the ray $[\gamma_\eta(N), \eta)$  and the $\delta$-ball around $\gamma_\eta(N)$ are disjoint from $D(\xi)$. Since $\eta \notin V_\cU(\xi)$, the $\xi$-family $\cU$ does not see $\gamma_\eta(N)$. Since $\cU'$ is %$((\lambda+1)\delta + \varepsilon)$-
refined in $\cU$ and since any geodesic between $\gamma_\eta(N)$ and a point of $B(\gamma_\eta(N), \delta)$ is contained in $B(\gamma_\eta(N), \delta)$, the $\xi$-family $\cU'$ does not see $B(\gamma_\eta(N), \delta)$, hence it does not see the following neighbourhood of $\eta$ 
$$U' := \{ y \in \overline{X} \mbox{ such that } \langle y, \eta \rangle_{v_0} > N \}.$$
Thus, $V_{U'}(\eta)$ is disjoint from $V_{\cU'}(\xi)$. 

Let $\xi' \in \partial_{Stab}G \setminus V_\cU(\xi)$ and let $\cU''$ be a $\xi'$-family that does not see $D(\xi)$. Since for every vertex $v$ of $D(\xi)$, we have $\xi' \notin \overline{U'_v}$, we can further assume that $\cU''$ and $\cU'$ are disjoint. Let us now prove that $V_{\cU'}(\xi)$ and $V_{\cU''}(\xi')$ are disjoint. Suppose by contradiction that this is not the case, and let $z \in \cC G$ be a point in that intersection and $x\in \overline{X}$ be a point of $p(z)$.  
First note that $x \notin D(\xi)$. Indeed, if that was the case, then since $\cU''$ does not see $D(\xi)$, we would have $x \in D(\xi) \cap D(\xi')$, which is absurd since  $\cU''$ and $\cU'$ are disjoint. Moreover, A geodesic from $v_0$ cannot leave $D(\xi)$ and $D(\xi')$ at the same vertex $v \in D(\xi) \cap D(\xi')$ for otherwise we would have $U_v \cap U_v' \neq \varnothing$, which is absurd. Thus a geodesic from $v_0$ to $x$ meets $D(\xi')$ after leaving $D(\xi)$. Since $\cU'$ is %$d_{max}$-
refined in $\cU$, it follows from the Refinement Lemma \ref{refinement}  that either there exists a vertex $v \in D(\xi) \cap D(\xi')$  such that $U_v \cap U_v' \neq \varnothing$, which is absurd by construction of $\cU''$, or we have $D(\xi')  \subset \mbox{Cone}_{\cU}(\xi)$, from which we deduce that   $\xi' \in V_\cU(\xi)$, a contradiction. 
\end{proof}

Note that Theorem \ref{basis} and Propositions \ref{T0}, \ref{regular}  immediately imply the following:

\begin{cor}[see {\cite[Theorem 7.12]{MartinBoundaries}}]
	The space $\cC G$ is metrizable. \qed
	\label{metrizable}
	\end{cor}

\subsection{Compactness}

We now prove the following: 

\begin{thm}[see {\cite[Theorem 7.13]{MartinBoundaries}}]
The space $\cC G$ is compact.
\label{compact}
\end{thm}

As $\cC G$ is metrisable by Corollary \ref{metrizable}, we just have to prove that $\cC G$ is sequentially compact. Let $(z_n)$ be a sequence in $\cC G$. As $\cV G$ is dense in $\cC G$, it is enough to consider a sequence $(z_n)$ of points of $\cV G$. For every $n$, set $x_n:= p(z_n) \in X$ and let $\gamma_n$ be a geodesic from $v_0$ to $x_n$, which  we write as a sequence $e_1^{(n)}, e_2^{(n)}, \ldots$. The proof splits in two cases.

\begin{lem}[see {\cite[Lemma 7.14]{MartinBoundaries}}]
Suppose that for every $k \geqslant 1$, the set $\{ e_k^{(n)}, n \geqslant 0\}$ is finite.
\begin{itemize}
\item If $d(v_0,x_n) \ra \infty$, then $(z_n)$ converges to a point of $\partial X$.
\item Otherwise, there exists a subsequence of $(z_n)$ converging to a point of $\cV G \cup \partial_{Stab} G$.
\end{itemize}
\label{compact1}
\end{lem}

\begin{proof}
Up to a subsequence, we can assume that there exists a sequence of edge $e_1, e_2, \ldots$ such that for every $k$, the sequence $(e_k^{(n)})$ is eventually constant at $e_k$. 
\begin{itemize}
\item The sequence $e_1, e_2, \ldots$ defines a $(\lambda,\varepsilon)$--quasi-geodesic $\gamma$ of $X$, since every finite segment $e_1, \ldots, e_k$ is eventually contained in $\gamma_n$ for $n$ large. Let $\eta$ be the associated point of $\partial X$. Since $\gamma_n$ and $\gamma$ share longer and longer initial segments, we have $\langle \gamma_n, \gamma\rangle_{v_0} \ra \infty$, hence $x_n$ converges to $\eta$ in $\overline{X}$. By definition of the topology of $\partial G$, this implies that $z_n$ converges to $\eta$ in $\partial G$. 
\item  Up to a subsequence, we can assume that $d(v_0, x_n)$ is constant and that $x_n$ is always a fixed vertex $v$ of $X$. Therefore, $z_n$ restrict to a sequence in $\cV G_v$. Either there exists a constant subsequence, in which case $z_n$ admits a subsequence converging to a point of $\cV G$, or we can find a subsequence of $(z_n)$ converging to a point $\xi \in \partial G_v$. In the latter case, Proposition \ref{inducedtopology} implies that this subsequence converges to $\xi$ in $\partial G$.\qedhere
\end{itemize}
\end{proof}

\begin{lem}[see {\cite[Lemma 7.15]{MartinBoundaries}}]
Suppose that there exists $k \geqslant 1$ such that the set $ \{ e_k^{(n)}, n \geqslant 0\}$ is infinite. Then there exists a subsequence of $(z_n)$ converging to a point of $\partial_{Stab} G$.
\label{compact2}
\end{lem}

\begin{proof}
Up to a subsequence, we can assume that geodesics $\gamma_n$ all start with edges $e_1, \ldots, e_{k-1}$ and the sequence of edges $e_k^{(n)}$ is injective. Let $v$ be their common vertex. Up to a subsequence, we can assume by cocompactness of the action that the sequence of subgroups $ G_{e_k^{(n)}} \subset \partial G_v$ is of the form $g_nG_eg_n^{-1}$ for some local group $G_e$ of $G(\Gamma)$ and $g_n \in G_v$, where the $g_n$ are in different cosets of $G_e$.  It follows from \cite[Theorem 1.8]{DahmaniCombination} that we can take a subsequence such that the sequence   $ \partial G_{e_k^{(n)}} \subset \partial G_v$ converges to an element $\xi \in \partial G_v$ in $\partial G_v$. Thus, for every $\xi$-family $\cU$, the point $x_n$ is eventually in $\mbox{Cone}_\cU(\xi)$ by the Geodesic Reattachment Lemma  \ref{geodesicreattachment}.
 This implies that $z_n$ is eventually in $V_\cU(\xi)$, hence $z_n$ converges to $\xi$.
\end{proof}

Following the same line of argument, we also get the following convergence criterion. 

\begin{cor}[see {\cite[Corollary 7.16]{MartinBoundaries}}]
Let $(K_n)$ be a sequence of subsets of $\cC G$.
\begin{itemize}
\item The sequence $(K_n)$ uniformly converges to a point $\eta \in \partial X$ if and only if the associated sequence of projections $p(K_n)$ uniformly converges to $\eta$ in $\overline{X}$.
\item Suppose that there exists a point $\xi$ of $\partial_{Stab} G$ such that, for $n$ large enough, every geodesic from $v_0$ to a point of $p(K_n)$ goes through $D(\xi)$. For every such $n$, every $z \in K_n$, choose a point $x \in \bar{p}(z)$ and a geodesic $\gamma_{v_0,x}$ and let $e_{n}(x)$ be the first edge touched by $\gamma_{v_0,x}$ after leaving $D(\xi)$. If there exists a vertex $v \in D(\xi)$ contained in each edge of the form $e_{n}(x)$ and such that for every neighbourhood $U$ of $\xi$ in $\cC G_v$, there exists an integer $N \geqslant 0$ such that for every $(n,x) \in \cup_{n \geqslant N} \{ n\} \times K_n$, we have $\cC G_{e_{n}(x)} \subset U$, then the sequence $(K_n)$ uniformly converges to $\xi$. \qed
\end{itemize}
\label{convergencecriterion}
\end{cor}

We now show that $G$ is a uniform convergence group on $\partial G$.  To that end, we follow closely the structure of \cite[Section 7]{MartinBoundaries}, and we refer to the corresponding result whenever possible.

\subsection{Continuity of the action}

We first show the following lemma. 

\begin{lem}[see {\cite[Lemma 6.18]{MartinBoundaries}}]
The topology of $\cC G$ (and $\partial G$) does not depend on the basepoint.
\label{basepoint}
\end{lem}

From now on, when dealing with neighbourhoods based at a given vertex, we may indicate that vertex as a superscript to avoid confusions.

\begin{proof} Let $v, v'$ two vertices of $X$. As the topology of $\overline{X}$ does not depend on the basepoint, we only consider the case of a point $\xi \in \partial_{Stab} G$. Let $\cU$ be a $\xi$-family for the topology based at $v$, and let $\cU'$ be a $\xi$-family for the topology based at $v'$ which is %$2D$-
refined in $\cU$. %(recall that $D$ was defined in ???). 
Let $x$ be a point of $\mbox{Cone}_{\cU'}^{v'}(\xi)$. By the Geodesic Reattachment Lemma \ref{geodesicreattachment}, any geodesic from $v$ to $x$ meets $D(\xi)$. By the Small Angle Property applied  to complexes $D(\xi)$, $\{x\}$ and to subsegments of geodesics from $v$ to $x$ and from $v'$ to $x$, it follows from the fact that $\cU'$  is 
 refined in $\cU$ that $x \in \mbox{Cone}_\cU^v(\xi)$. Since $\cU'$ is contained in $\cU$, it follows that $V_{\cU'}^{v'}(\xi) \subset V_\cU^v(\xi)$.
\end{proof}

The action of $G$ on $EG$ extends to $\partial G$ as follows. The complex $X$ being hyperbolic, the action of $G$ on $X$ by isometries extends to $\partial X$. Furthermore, we described in Section \ref{s:constrbdry} an action of $G$ on $\partial_{Stab} G$. This defines a $G$-action on $\partial G$.

\begin{cor}
The action of $G$ on $\cC G$ and $\partial G$ is continuous. 
\label{continuousaction}
\end{cor}

\begin{proof}
Let $g$ be an element of $G$ and $z$ a point of $\cC G$. The element $g$ sends a basis of neighbourhoods of $z$ centred at a vertex $v$ to a basis of neighbourhoods of $\xi$ centred at $gv$. Since the topology of $\partial G$ does not depend on the basepoint by Lemma \ref{basepoint}, the result follows.
\end{proof}
\subsection{Convergence group action}

The remainder of the proof, namely that $G$ acts as a uniform convergence group on the boundary $\partial G$, is completely analogous to the proof presented in \cite[Section 9]{MartinBoundaries}. Indeed, the proof relied on the CAT(0) version of the Geodesic Reattachment Lemma \ref{geodesicreattachment} and Refinement Lemma \ref{refinement}. As we introduced combinatorial analogues of these tools in Section \ref{sec:toolbox}, the proof carries over to this combinatorial setting without any essential change. The only instance where an additional notion was being used in the statement is in the following: 

\begin{lem}[see {\cite[Lemma 9.16]{MartinBoundaries}}] Let $(g_n)$ be an injective sequence of elements of $G$, and suppose that for some (hence every) vertex $v_0$ of $X$, we have $d(v_0, g_nv_0) \ra \infty$. Since $(EG \cup \partial G, \partial G)$ is an $E\cZ$-structure, we assume that there exist $\xi_+, \xi_- \in \partial G$  such that for every compact $K \subset EG$, we have $g_n K \ra \xi_+ $ and  $g_n^{-1} K \ra \xi_+ $. Then there exists  a subsequence $(g_{\varphi(n)})$ such that for every compact subset $K$ of $\partial G \setminus \left\{\xi_-\right\}$, the sequence of translates $g_{\varphi(n)}K$ uniformly converges to $\xi_+$.
\end{lem}

In this statement, the use of $E\cZ$-structures was convenient but unnecessary. Indeed, let $x_0$ be a point of $\cV G$. %Since $(\overline{EG}, \partial G)$ is an $E\cZ$-structure for $G$ by Theorem \ref{EZstructure}, 
By compactness of $\cC G$ and up to taking a subsequence, we can assume that there exist $\xi_+, \xi_- \in \partial G$ such that %for every compact subset $K \subset EG$, 
%we have
$g_n x_0 \ra \xi_+$ and $g_n^{-1} x_0 \ra \xi_-$. For every point $x$ of $\cV G$, we also have $g_n x \ra \xi_+$ and $g_n^{-1} x \ra \xi_-$. This is clear if $\xi_{\pm}$ is in $\bd X$ since $g_n^{\pm}v_0$ and $g_n^{\pm}v$ stay at bounded distance in $X$. If $\xi_{\pm}$ is in $\partial_{Stab} G$, this follows from the Refinement Property \ref{refinement} since the interval $I(g_n^{\pm}v_0,g_n^{\pm}v)$ is disjoint from $D(\xi_{\pm})$ for $n$ big enough, as $d(v_0, g_n^{\pm} v_0) \ra \infty$. Thus, in our setting, we are led to prove the following analogous statement: 

\begin{lem}Let $(g_n)$ be an injective sequence of elements of $G$, and suppose that for some (hence every) vertex $v_0$ of $X$, $d(v_0, g_n v_0) \ra \infty$.  %Since $(\overline{EG}, \partial G)$ is an $E\cZ$-structure for $G$ by Theorem \ref{EZstructure}, 
	Up to taking a subsequence, we assume that there exist $\xi_+, \xi_- \in \partial G$ such that for every $x \in \cV G$, 
	we have
	$g_n x \ra \xi_+$ and $g_n^{-1} x \ra \xi_-$. 
	
	Then there exists a subsequence $(g_{\varphi(n)})$ such that for every compact subset $K$ of $\partial G \setminus \left\{\xi_-\right\}$, the sequence of translates $g_{\varphi(n)}K$ uniformly converges to $\xi_+$. 
	\label{convergencegroup3}
\end{lem}

The proof of this statement is then completely analogous to that of \cite[Lemma 9.16]{MartinBoundaries}, and relies on the convergence criterion \ref{convergencecriterion}.

\end{appendices}